\documentclass[a4paper,12pt]{article}

\usepackage{latexsym}
\usepackage{amssymb}
\usepackage{theorem}
\usepackage{amsmath}
\usepackage{amscd}
\usepackage{graphicx}

\newtheorem{theorem}{Theorem}[section]
\newtheorem{corollary}[theorem]{Corollary}

\newtheorem{example}[theorem]{Example}
\newtheorem{proposition}[theorem]{Proposition}
\newtheorem{remark}[theorem]{Remark}
\newtheorem{definition}[theorem]{Definition}

\newcommand{\demo}{\par\noindent{\it Proof. \/}\ }
\newcommand{\enD}{\hfill $\Box$\vspace{3truemm} \par}
\newcommand{\R}{\mathbb{R}}

\newcommand{\bn}{\mbox{\boldmath $n$}}
\newcommand{\bt}{\mbox{\boldmath $t$}}
\newcommand{\ba}{\mbox{\boldmath $a$}}
\newcommand{\bb}{\mbox{\boldmath $b$}}

\newcommand{\bx}{\mbox{\boldmath $x$}}
\newcommand{\A}{\mathcal{A}}

\begin{document}

\title{On vertices of frontals in the Euclidean plane}

\author{Nozomi Nakatsuyama and Masatomo Takahashi}

\date{\today}

\maketitle

\begin{abstract}
We investigate vertices for plane curves with singular points. 
As plane curves with singular points, we consider Legendre curves (respectively, Legendre immersions) in the unit tangent bundle over the Euclidean plane and frontals (respectively, fronts) in the Euclidean plane. 
We define a vertex using evolutes of frontals. 
After that we define a vertex of a frontal in the general case. 
It is also known that the four vertex theorem does not hold for simple closed fronts. 
We give conditions under which a frontal has a vertex and the four vertex theorem holds for closed frontals. 
We also give examples and counter examples of the four vertex theorem.
\end{abstract}

\renewcommand{\thefootnote}{\fnsymbol{footnote}}
\footnote[0]{2020 Mathematics Subject classification: 58K05, 53A04, 57R45}
\footnote[0]{Key Words and Phrases. vertex, Legendre curve, frontal, inflection point, evolute, singular point}

\section{Introduction}

The four vertex theorem is a classical result in differential geometry. 
It states that there are four vertices for simple closed regular convex plane curves (cf.  \cite{Mukhopadhyaya}). 
In general, there are four vertices for simple closed regular plane curves (cf. \cite{Kneser}). 
There are counter examples for non-simple closed regular plane curves. 
See also \cite{Dias-Tari, Gounai-Umehara,Salarinoghabi-Tari,Tabachnikov,Thorbergsson-Umehara,Umehara,Umehara-Yamada}. 
We investigate vertices for plane curves with singular points. 
As smooth curves with singular points, we consider Legendre curves (respectively, Legendre immersions) in the unit tangent bundle over the Euclidean plane and frontals (respectively, fronts) in the Euclidean plane. 
The frontals form a considerably large class of singular curves, so
the results in this paper are not too specific.
For basic results on singularity theory see \cite{Arnold1,Arnold2,Bruce-Giblin,Ishikawa-book,Izumiya-book}. 
It is also known that the four vertex theorem does not hold for simple closed fronts (see, Example \ref{nephroid}). 
In \cite{Fukunaga-Takahashi-2014}, we gave conditions for a front to have at least four vertices. 
In this paper, we give conditions under which a frontal has a vertex and the four vertex theorem holds for closed frontals.  
\par
In \S 2, we review the theories of regular curves, Legendre curves on the unit tangent bundle over the Euclidean plane and the evolutes of fronts. 
In \S 3, using evolutes of frontals, we define the vertex as a singular point of the evolute and investigate the four vertex theorem of frontals. 
We give conditions that frontal has at least four vertices using the curvature of the Legendre curve and the type of singular points of the Legendre curve. 
Moreover, for a closed Legendre curve if inflection and singular points of the frontal are isolated points, and the frontal is a simple convex frontal, then the frontal has at least four vertices (Theorem \ref{four-vertex-theorem-frontal}). 
In \S 4 we define a vertex of a frontal in the general case. 
We investigate properties of vertices of frontals. 
We show that vertices of frontals are affine invariant (Propositions \ref{affine1} and \ref{evolute-affine.frontal}).
For a closed Legendre immersion, if the front is a simple convex with at least two singular and two inflection points, then the front has at least four vertices (Theorem \ref{four-vertex-theorem-front-inflection}).

We shall assume throughout the whole paper that all maps and manifolds are $C^{\infty}$ unless the contrary is explicitly stated.

\bigskip
\noindent
{\bf Acknowledgement}. 
The second author was supported by JSPS KAKENHI Grant Number No. 20K03573.

\section{Preliminaries}

We quickly review the theories of regular curves, Legendre curves on the unit tangent bundle over $\R^2$ and the evolutes of fronts. 
\par
Let $I$ be an interval or $\R$, and let $\R^2$ be the Euclidean plane with the inner product $\ba \cdot \bb=a_1b_1+a_2b_2$, where $\ba=(a_1,a_2)$ and $\bb=(b_1,b_2) \in \R^2$. We denote $|\ba|=\sqrt{\ba\cdot\ba}$.

\subsection{Regular plane curves}

Let $\gamma: I \to \R^2$ be a regular plane curve, that is, $\dot{\gamma}(t)=(d\gamma/dt)(t) \not=0$ for all $t \in I$. 
We have the unit tangent vector $\bt(t)=\dot{\gamma}(t)/|\dot\gamma(t)|$ and the unit normal vector $\bn(t)=J(\bt(t))$, where $J$ is the anti-clockwise rotation by $\pi/2$ on $\R^2$. 
Then we have the Frenet formula
$$
\left(
\begin{array}{c}
\dot{\bt}(t)\\
\dot{\bn}(t)
\end{array}
\right)
=
\left(
\begin{array}{cc}
0 & |\dot\gamma(t)|\kappa(t)\\
-|\dot\gamma(t)|\kappa(t) & 0
\end{array}
\right)
\left(
\begin{array}{c}
\bt(t)\\
\bn(t)
\end{array}
\right),
$$
where $\kappa:I \to \R$ is the curvature.
Note that the curvature $\kappa(t)$ is independent of the choice of a parametrization up to sign.

We say that $t_0 \in I$ is an {\it inflection point} of $\gamma$ if $\kappa(t_0)=0$ and 
a {\it vertex} of $\gamma$ if $\dot{\kappa}(t_0)=0$.
Then the classical four vertex theorem is as follows.

\begin{theorem}[The four vertex theorem of regular curves, \cite{Kneser, Mukhopadhyaya}]
If $\gamma:I \to \R^2$ is a simple convex regular curve, then $\gamma$ has at least four vertices.
In general, if $\gamma:I \to \R^2$ is a simple closed regular curve, then $\gamma$ has at least four vertices.
\end{theorem}

Suppose that $\gamma$ does not have inflection points. 
We define an {\it evolute} $Ev(\gamma):I \to \R^2$ of $\gamma$ by 
$
Ev(\gamma)(t)=\gamma(t)+(1/\kappa(t))\bn(t).
$
Then $t_0$ is a singular point of $Ev(\gamma)$ if and only if $t_0$ is a vertex of $\gamma$.
It follows that $Ev(\gamma)$ has at least four singular points for simple closed regular curves $\gamma$ without inflection points (see Example \ref{ellipse}). 
We investigate this phenomena for curves with singular points by using evolutes. 
\begin{remark}{\rm 
Suppose that $t_0$ is a regular point and not an inflection point of $\gamma$.
It is well-known that if $t_0$ is a vertex of $\gamma$, then $\gamma$ and the osculating circle have at least third order contact at $t_0$ (cf. \cite{doCarmo,Gibson,Umehara-Yamada}). 
}
\end{remark}
\subsection{Legendre curves in the unit tangent bundle over the Euclidean plane}

Let $(\gamma,\nu): I \to \R^2 \times S^1$ be a smooth mapping, where  $S^1$ is the unit circle. 
We say that $(\gamma,\nu):I \to \R^2 \times S^1$ is {\it a Legendre curve} if 
$(\gamma,\nu)^* \theta=0$, where $\theta$ is the canonical contact $1$-form on the unit tangent bundle $T_1 \R^2=\R^2 \times S^1$ (cf. \cite{Arnold1, Arnold2}). 
This condition is equivalent to $\dot{\gamma}(t) \cdot \nu(t)=0$ for all $t \in I$. 
Moreover, if $(\gamma,\nu)$ is an immersion, we call $(\gamma,\nu)$ {\it a Legendre immersion}.
We say that $\gamma: I \to \R^2$ is {\it a frontal} (respectively, {\it a front} or {\it a wave front}) if there exists a smooth mapping $\nu:I \to S^1$ such that 
$(\gamma,\nu)$ is a Legendre curve (respectively, a Legendre immersion).

Let $(\gamma,\nu):I \to \R^2 \times S^1$ be a Legendre curve.
We put on $\mu (t)=J(\nu (t))$ and call the pair $\{\nu(t), \mu(t) \}$ a {\it moving frame along the frontal $\gamma(t)$} in $\R^2$. 
Then we have the Frenet type formula of the Legendre curve which is given by
\begin{eqnarray*}\label{Frenet.frontal}
\left(
\begin{array}{c}
\dot{\nu}(t)\\
\dot{\mu}(t)
\end{array}
\right)
=
\left(
\begin{array}{cc}
0 & \ell(t)\\
-\ell(t) & 0
\end{array}
\right)
\left(
\begin{array}{c}
\nu(t)\\
\mu(t)
\end{array}
\right), \ 
\dot\gamma(t)=\beta(t) \mu(t),
\end{eqnarray*}
where $\ell(t)=\dot\nu(t) \cdot \mu(t), \beta(t)=\dot{\gamma}(t) \cdot \mu(t)$. 
We call the pair $(\ell,\beta)$ {\it the (Legendre) curvature of the Legendre curve} (cf. \cite{Fukunaga-Takahashi-2013}).
Then $(\gamma,\nu):I \to \R^2 \times S^1$ is a Legendre immersion if and only if  $(\ell(t),\beta(t)) \not=(0,0)$ for all $t \in I$. We have the existence and the uniqueness for Legendre curves in the unit tangent bundle like as regular plane curves, see \cite{Fukunaga-Takahashi-2013}.

\begin{example}[Regular curves]\label{regular-curves}{\rm
Let $\gamma:I \to \R^2$ be a regular curve. 
If we consider $\bn$ as a unit normal vector, then $(\gamma,\bn):I \to \R^2 \times S^1$ is a Legendre immersion with the curvature $(\ell,\beta)=(|\dot{\gamma}|\kappa,-|\dot{\gamma}|)$. 
It follows that $\kappa(t)=-\ell(t)/\beta(t)$. 
If we consider $-\bn$ as a unit normal vector, then $(\gamma,-\bn):I \to \R^2 \times S^1$ is also a Legendre immersion with the curvature $(\ell,\beta)=(|\dot{\gamma}|\kappa,|\dot{\gamma}|)$. 
It follows that $\kappa(t)=\ell(t)/\beta(t)$ (cf. Proposition \ref{nu-change}). 
}
\end{example}
\begin{example}[Type $(n,m)$]\label{nm-type}{\rm 
Let $n, k \in \mathbb{N}$ and $m=n+k$. 
We consider a smooth map germ $(\gamma,\nu):(\R,0) \to \R^2 \times S^1$,
\begin{align*}
\gamma(t) =(\pm t^n,t^mf(t)), \
\nu(t) =\frac{(-mt^kf(t)-t^{k+1}\dot{f}(t),\pm n)}{\sqrt{(mt^kf(t)+t^{k+1}\dot{f}(t))^2+n^2}},
\end{align*}
where $f:(\R,0) \to \R$ is a smooth function germ with $f(0) \not=0$. 
Note that $0$ is a singular point of $\gamma$ when $n>1$. 
Then $(\gamma,\nu)$ is a Legendre curve with curvature 
\begin{align*}
\ell(t)&=\pm \frac{nt^{k-1}(mkf(t)+(m+k+1)t\dot{f}(t)+t^{2}\ddot{f}(t))}{(mt^kf(t)+t^{k+1}\dot{f}(t))^2+n^2},\\ \beta(t)&=-t^{n-1}\sqrt{(mt^kf(t)+t^{k+1}\dot{f}(t))^2+n^2}.
\end{align*}
We say $\gamma$ is of type $(n,m)$. 
If $n=1$, then $\gamma$ is regular (a front) and $(\gamma,\nu)$ is a Legendre immersion.
If $k=1$, then $\gamma$ is a front and $(\gamma,\nu)$ is a Legendre immersion.
Otherwise, that is, if $n,k>1$, then $\gamma$ is a frontal and $(\gamma,\nu)$ is a Legendre curve.
}
\end{example}

\begin{proposition}[\cite{Fukunaga-Takahashi-2013}]\label{Legendre.function}
Let $(\gamma,\nu):I \to \R^2 \times S^1$ be a Legendre curve with curvature $(\ell,\beta)$ and $t:\widetilde{I} \to I$ be a smooth function, where $\widetilde{I}$ is an interval. Then $(\gamma \circ t, \nu \circ t): \widetilde{I} \to \R^2 \times S^1$ is also a Legendre curve with the curvature $((\ell \circ t)t', (\beta \circ t)t')$.
\end{proposition}
We say that $t:\widetilde{I} \to I$ is a parameter change if $t'(u) \not=0$ for all $u \in \widetilde{I}$. 
Two map germs $f,g\colon(\R^n,0)\to (\R^p,0)$ are said to be $\A$-equivalent
if there exist diffeomorphism germs 
$\phi\colon (\R^n,0)\to(\R^n,0)$
and
$\Phi \colon (\R^p,0)\to(\R^p,0)$
such that 
$f=\Phi \circ g\circ \phi$.
When $\Phi$ is the identity in the above,
$f$ and $g$ are said to be $\mathcal{R}$-equivalent. 
See \cite{Arnold1, Arnold2, Bruce-Giblin, Izumiya-book} for basic idea of singularity theory.
The $\mathcal{A}$ or $\mathcal{R}$-equivalent class of mappings or functions
are worth studying from the viewpoint of differential topology. 
A one-variable function $f:(I,t_0) \to (\R,0)$ has type $A_k$ at $t_0 \in I$ if $f^{(i)}(t_0)=0$ for $i=1,\dots,k$ and $f^{(k+1)}(t_0) \not=0$. 
Then $f$ is $\mathcal{R}$-equivalent to $\pm t^{k+1}$ (cf. \cite{Bruce-Giblin}).

\begin{proposition}\label{nm-type-parametrization}
Let $\gamma:(I,t_0) \to (\R^2,0)$ be a smooth map germ with $\gamma(t)=(x(t),y(t))$ and $2 \le n < m$. 
Suppose that $x$ has type $A_{n-1}$ and $y$ has type $A_{m-1}$ at $t_0$. 
Then $\gamma$ is $\mathcal{R}$-equivalent to of type $(n,m)$.
\end{proposition}

We consider how to change the curvature of Legendre curves by a diffeomorphism of the target.
\begin{proposition}\label{diffeomorphism-target}
Let $(\gamma,\nu):I \to \R^2 \times S^1$ be a Legendre curve with curvature $(\ell,\beta)$ and 
$\Phi:\R^2 \to \R^2, \Phi(x,y)=(\phi_1(x,y),\phi_2(x,y))$ be a diffeomorphism. 
We denote $\gamma(t)=(x(t),y(t)), \nu(t)=(a(t),b(t))$ and $\widetilde{\gamma}=\Phi \circ \gamma$.
Then $(\widetilde{\gamma}, \widetilde{\nu}):I \to \R^2 \times S^1$ is a Legendre curve, where $\widetilde{\nu}=\overline{\nu}/|\overline \nu|$,
$$
\overline{\nu}(t)=(\phi_{2y}(\gamma(t))a(t)-\phi_{2x}(\gamma(t))b(t),-\phi_{1y}(\gamma(t))a(t)+\phi_{1x}(\gamma(t))b(t)).
$$
The curvature $(\widetilde{\ell},\widetilde{\beta})$ is given by 
\begin{align*} 
\widetilde{\ell}(t) &= \frac{1}{|\overline{\nu}(t)|^2}\Bigl( \bigl((\phi_{2xx}(\gamma(t))b^2(t)-2\phi_{2xy}(\gamma(t))a(t)b(t)+\phi_{2yy}(\gamma(t))a^2(t))\\
& \qquad (-\phi_{1x}(\gamma(t))b(t)+\phi_{1y}(\gamma(t))a(t))\\
& \quad -(\phi_{1xx}(\gamma(t))b^2(t)-2\phi_{1xy}(\gamma(t))a(t)b(t)+\phi_{1yy}(\gamma(t))a^2(t))\\
& \qquad (-\phi_{2x}(\gamma(t))b(t)+\phi_{2y}(\gamma(t))a(t)) \bigr)\beta(t)\\
& \quad +(\phi_{1x}(\gamma(t))\phi_{2y}(\gamma(t))-\phi_{2x}(\gamma(t))\phi_{1y}(\gamma(t)))\ell(t) \Bigr),\\
\widetilde{\beta}(t)&=|\overline{\nu}(t)|\beta(t).
\end{align*}
\end{proposition}
\demo
Since $\dot{\widetilde{\gamma}}(t)=(\phi_{1x}(\gamma(t))\dot{x}(t)+\phi_{1y}(\gamma(t))\dot{y}(t),  \phi_{2x}(\gamma(t))\dot{x}(t)+\phi_{2y}(\gamma(t))\dot{y}(t))$, 
we have $\dot{\widetilde{\gamma}}(t) \cdot \overline{\nu}(t)=0$ for all $t \in I$, where 
$$
\overline{\nu}(t)=(\phi_{2y}(\gamma(t))a(t)-\phi_{2x}(\gamma(t))b(t),-\phi_{1y}(\gamma(t))a(t)+\phi_{1x}(\gamma(t))b(t)).
$$ 
Moreover, since $\Phi$ is a diffeomorphism, $|\overline{\nu}(t)| \not=0$ for all $t \in I$. 
It follows that $(\widetilde{\gamma},\widetilde{\nu}):I \to \R^2 \times S^1$ is a Legendre curve, where $\widetilde{\nu}=\overline{\nu}/|\overline \nu|$.
By a direct calculation, we have the curvature $(\widetilde{\ell},\widetilde{\beta})$ of $(\widetilde{\gamma},\widetilde{\nu})$.
\enD
As special cases, we consider affine transformations and the reflection of the target of $\gamma$.
\begin{corollary}\label{change-target}
Let $(\gamma,\nu):I \to \R^2 \times S^1$ be a Legendre curve with curvature $(\ell,\beta)$ and $\Phi:\R^2 \to \R^2$ be a diffeomorphism. 
We denote $\nu=(a,b)$ and $\widetilde{\gamma}=\Phi \circ \gamma$.
\par
$(1)$ Suppose that the diffeomorphism $\Phi:\R^2 \to \R^2$ is given by  $\Phi(x,y)=(a_{11}x+a_{12}y,a_{21}x+a_{22}y)$, where $a_{11}a_{22}-a_{12}a_{21} \not=0$ and $a_{11},a_{12},a_{21},a_{22} \in \R$. 
Then $(\widetilde{\gamma},\widetilde{\nu}):I \to \R^2 \times S^1$ is a Legendre curve with the curvature $(\widetilde{\ell},\widetilde{\beta})=\left((a_{11}a_{22}-a_{12}a_{21})\ell/|\overline{\nu}|^2,|\overline{\nu}|\beta\right)$, where $\widetilde{\nu}=\overline{\nu}/|\overline{\nu}|$ and 
$\overline{\nu}=(a_{22}a-a_{21}b,-a_{12}a+a_{11}b)$.
\par
$(2)$ Suppose that the diffeomorphism $\Phi: \R^2 \to \R^2$ is given by $\Phi(x,y)=(y,x)$.  
Then $(\widetilde{\gamma},\widetilde{\nu}):I \to \R^2 \times S^1$ is a Legendre curve with the curvature $(\widetilde{\ell},\widetilde{\beta})=\left(-\ell,\beta\right)$, where $\widetilde{\nu}=(-b,-a)$.
\end{corollary}
By a direct calculation, we have the following.
\begin{proposition}\label{nu-change}
Let $(\gamma,\nu):I \to \R^2 \times S^1$ be a Legendre curve with curvature $(\ell,\beta)$. 
Then $(\gamma,-\nu):I \to \R^2 \times S^1$ is also a Legendre curve with the curvature $(\ell,-\beta)$.
\end{proposition}

Let $(\gamma,\nu):I \to \R^2 \times S^1$ be a Legendre curve with curvature $(\ell,\beta)$.
We say that a point $t_0 \in I$ is an {\it inflection point} of  the frontal $\gamma$ (or, the Legendre curve $(\gamma,\nu)$) if $\ell(t_0)=0$.
Remark that the definition of the inflection point of the frontal is a generalisation of the definition of the inflection point of a regular curve by Example \ref{regular-curves} (cf. \cite{Fukunaga-Takahashi-2013}).
\par
We define a {\it parallel curve} $\gamma^\lambda: I \to \R^2$ of the frontal $\gamma$ (or, Legendre curve $(\gamma,\nu)$) by 
$
\gamma^\lambda(t)=\gamma(t)+\lambda \nu(t),
$
where $\lambda \in \R$.
\begin{proposition}[\cite{Fukunaga-Takahashi-2015}]\label{parallel.frontal}
Let $(\gamma, \nu) : I \to \mathbb{R}^2 \times S^1$ be a Legendre curve with curvature $(\ell, \beta)$. 
The parallel curve $\gamma^\lambda:I \to \R^2$ is also a frontal for any $\lambda \in \R$. 
More precisely, $(\gamma^\lambda,\nu):I \to \R^2 \times S^1$ is a Legendre curve with the curvature
$\left(\ell, \beta+\lambda \ell \right).$
\end{proposition}

For $n \in \mathbb{N} \cup \{0\}$, we say that a Legendre curve $(\gamma, \nu) : [a,b] \rightarrow \mathbb{R}^2 \times S^1$ is $C^{n}$-\emph{closed} if $(\gamma^{(k)}(a), \nu^{(k)}(a))=(\gamma^{(k)}(b), \nu^{(k)}(b))$ for all $ k \in \{0,\cdots,n \}$, where $\gamma^{(k)}(a)$, $\nu^{(k)}(a)$,  $\gamma^{(k)}(b)$ and $\nu^{(k)}(b)$ mean one-sided $k$-th differential.
In this paper, we say that $(\gamma, \nu)$ is a \emph{closed} Legendre curve, if the
curve is at least $C^{1}$-closed (cf. \cite{Fukunaga-Takahashi-2016}). When $a$ and $b$ are singular points of $\gamma$, we treat these singular points as one singular point. 
Moreover, a frontal $\gamma : [a,b] \rightarrow \mathbb{R}^2$ is {\em simple closed} if for $t_1 < t_2$, we have $\gamma(t_1)=\gamma(t_2)$ if and only if $t_1 = a$ and $t_2 = b$. 
We define a convex frontal in the Euclidean plane. 
From now on, $I$ is a closed interval. Let $(\gamma, \nu) : I \rightarrow \mathbb{R}^2 \times S^1$ be a Legendre curve. 
We denote the tangent line at $t$ of $\gamma$ by $L_{t}$, that is, $L_{t} =\{ \lambda \mu(t)+\gamma(t) \  | \  \lambda \in \mathbb{R} \}$. 
Any tangent line $L_{t}$ divides $\mathbb{R}^2$ into two half-planes $H_+$ and $H_-$ such
that $H_+ \cup H_- = \mathbb{R}^2$ and $H_+ \cap H_- = L_{t}$. 
By using $\nu$, the half-planes $H_+$ and $H_-$ are presented by 
$H_+ = \{ \bx \in \mathbb{R}^2  \ | \ ( \bx - \gamma(t)) \cdot \nu(t) \ge 0 \}$  and 
$H_- = \{ \bx \in \mathbb{R}^2  \ | \ ( \bx - \gamma(t)) \cdot \nu(t) \le 0 \}$.
For a Legendre curve $(\gamma, \nu) : I \rightarrow \mathbb{R}^2 \times S^1$, 
we say that $(\gamma, \nu)$ is a {\em convex Legendre curve} (or, $\gamma$ is a {\em convex frontal}) if $\gamma(I) \subset H_+$ for any tangent line of $\gamma$ or $\gamma(I) \subset H_-$ for any tangent line of $\gamma$.
Note that if $\gamma$ is a regular curve, then $\mu(t)$ is equal to the unit
tangent vector of $\gamma$ at $\gamma(t)$ up to sign. Therefore, $\gamma$ is a convex curve as a frontal if and only if $\gamma$ is a convex curve as the usual mean when $\gamma$ is a regular curve (cf. \cite{Gray}). 
\begin{theorem}[\cite{Fukunaga-Takahashi-2016}]\label{Convex-frontal}
Let $(\gamma, \nu) : I \rightarrow \mathbb{R}^2 \times S^1$ be a closed Legendre curve with  curvature $(\ell,\beta)$ which the frontal $\gamma$ is simple closed. Suppose that zeros of $\ell$ and of $\beta$ are isolated points.
Then the frontal $\gamma$ is convex if and only if the curvature satisfy one of the following condition:
\par
${\rm (i)}$ Both of $\ell(t)$ and $\beta(t)$ are always non-negative,
\par  
${\rm (ii)}$ $\ell(t)$ is always non-negative and $\beta(t)$ is always non-positive,
\par
${\rm (iii)}$ Both of $\ell(t)$ and $\beta(t)$ are always non-positive,
\par
${\rm (iv)}$ $\ell(t)$ is always non-positive and $\beta(t)$ is always non-negative.  
\end{theorem}

\subsection{Evolutes of fronts}

In order to define a vertex of a front, we consider an evolute of the front, in detail see \cite{Fukunaga-Takahashi-2014}. 
Let $(\gamma,\nu):I \to \R^2 \times S^1$ be a Legendre curve with curvature  $(\ell,\beta)$.
If $(\gamma,\nu)$ does not have inflection points, namely, $\ell(t) \not=0$ for all $t \in I$, then $(\gamma,\nu)$ is a Legendre immersion. 
In this subsection, we assume that $(\gamma,\nu)$ does not have inflection points.
\begin{definition}[\cite{Fukunaga-Takahashi-2014}]\label{evo-front.def}{\rm 
The {\it evolute} $\mathcal{E}v(\gamma):I \to \R^2$ of the front $\gamma$ is given by 
\begin{eqnarray*}\label{evo-front}
\mathcal{E}v(\gamma)(t)=\gamma(t)-\frac{\beta(t)}{\ell(t)}\nu(t).
\end{eqnarray*}
}
\end{definition}

For a Legendre immersion $(\gamma,\nu)$ with curvature $(\ell,\beta)$, $t_0$ is {\it a vertex of the front $\gamma$} (or {\it a Legendre immersion $(\gamma,\nu)$})
if $(d/dt)\mathcal{E}v(t_0)=0$, namely, $(d/dt)(\beta/\ell)(t_0) =0$. 
Note that if $t_0$ is a regular point of $\gamma$, the definition of the vertex coincides with usual vertex for regular curves without inflection points by Example \ref{regular-curves}. 
\par
By the definition of parallel curves, $(\gamma^\lambda,\nu)$ is also a Legendre immersion without inflection points. 
Moreover, the evolute $\mathcal{E}v(\gamma^\lambda)$ coincides with the evolute $\mathcal{E}v(\gamma)$.

\begin{proposition}[\cite{Fukunaga-Takahashi-2014}]\label{Four-vertex-theorem-front}
Let $(\gamma,\nu):[a,b] \to \R^2 \times S^1$ be a closed Legendre immersion without inflection points.
\par
$(1)$ If $\gamma$ has at least two singular points which degenerate more
 than $3/2$ cusp, then $\gamma$ has at least four vertices.
\par
$(2)$ If $\gamma$ has at least four singular points, then $\gamma$ has 
at least four vertices.  
\end{proposition}
\begin{example}[Parallel curves of an ellipse]\label{ellipse}{\rm
Let $(\gamma,\nu):[0,2\pi] \to \R^2 \times S^1$ be 
\begin{align*}
\gamma(t)=(a \cos t,b \sin t), \ \nu(t)=\frac{1}{\sqrt{a^2 \sin^2 t+b^2 \cos^2 t}}(-b \cos t, -a \sin t).
\end{align*}
Here, we consider $a,b \in \R$ with $a>b>0$. 
By a direct calculation, $(\gamma,\nu)$ is a closed Legendre immersion with the curvature 
\begin{align*}
\ell(t)=\frac{ab}{a^2 \sin^2 t+b^2 \cos^2 t}, \ \beta(t)=-\sqrt{a^2 \sin^2 t+b^2 \cos^2 t}. 
\end{align*}
Note that $\gamma$ is a regular curve and a simple convex closed curve. 
Therefore, there exist at least four vertices of $\gamma$. 
In fact, we see that $t=0, {\pi}/{2}, \pi, 3\pi/2$ are vertices of $\gamma$.
The evolute of the ellipse $\gamma$ is given by 
\begin{align*}
\mathcal{E}v(\gamma)(t)=\left(\frac{a^2-b^2}{a}\cos^3 t,-\frac{a^2-b^2}{b}\sin^3 t\right).
\end{align*}
Moreover, the parallel curve for each $\lambda\in\R$ is given by 
\begin{align*}
\gamma^\lambda(t)=\left(\cos t\left(a-\frac{b\lambda}{\sqrt{a^2 \sin^2 t+b^2 \cos^2 t}}\right), \sin t\left(b-\frac{a\lambda}{\sqrt{a^2 \sin^2 t+b^2 \cos^2 t}}\right)\right). 
\end{align*}
If we take $\lambda=b^2/a$, then $\gamma^\lambda$ at $t=0$ and $t=\pi$ are diffeomorphic to the 4/3 cusps (cf. \cite{Bruce-Gaffney,Fukunaga-Takahashi-2014}). 
By Proposition \ref{Four-vertex-theorem-front} $(2)$, $\gamma^\lambda$ has (at least) four vertices.  
}
\end{example}
\begin{example}[Nephroid]\label{nephroid}{\rm
Let $(\gamma,\nu):[0,2\pi] \to \R^2 \times S^1$ be 
\begin{align*}
\gamma(t) =(3\cos t-\cos 3t,3\sin t-\sin 3t), \
\nu(t) =(-\sin 2t,\cos 2t).
\end{align*}
By a direct calculation, $(\gamma,\nu)$ is a closed Legendre immersion with the curvature $\ell(t)=2, \beta(t)=-6\sin t$. 
We call $\gamma$ the nephroid. 
Since $(d/dt)(\beta(t)/\ell(t))=-3\cos t$, $t={\pi}/{2}$ and ${3\pi}/{2}$ are vertices. 
Therefore, $\gamma$ is a simple closed front, but has only two vertices.
Moreover, we can show that $t=0$ and $\pi$ are $3/2$-cusps of $\gamma$ (cf. \cite{Bruce-Gaffney,Fukunaga-Takahashi-2014}).  
Furthermore, by Theorem \ref{Convex-frontal}, $\gamma$ is not convex. 
The evolute of the nephroid is given by
$\mathcal{E}v(\gamma)(t) =\left(2 \cos^3 t, 3 \sin t-2\sin^3 t \right).$
}
\end{example} 
\section{Evolutes of frontals and vertices}

In order to define a vertex of a frontal, we consider an evolute of the frontal, in detail see \cite{Fukunaga-Takahashi-2015}.
For more properties of frontal see in \cite{Ishikawa-2018,Ishikawa-2020}.
Let $(\gamma,\nu):I \to \R^2 \times S^1$ be a Legendre curve with curvature  $(\ell,\beta)$.

\begin{definition}[\cite{Fukunaga-Takahashi-2015}]\label{evolute-frontal.def}{\rm
The {\it evolute $\mathcal{E}v(\gamma):I \to \R^2$ of the frontal $\gamma$} is given by 
\begin{eqnarray*}\label{evo-frontal}
\mathcal{E}v(\gamma)(t)=\gamma(t)-\alpha(t)\nu(t),
\end{eqnarray*}
if there exists a unique smooth function $\alpha:I \to \R$ such that $\beta(t)=\alpha(t)\ell(t)$.
In this case, we say that {\it the evolute $\mathcal{E}v(\gamma)$ exists}.
}
\end{definition}

If $t_0$ is an inflection point of $(\gamma, \nu)$ and the evolute  $\mathcal{E}v(\gamma)$ exists, then the inflection point must be a singular point of $\gamma$. 
It follows that $t_0$ is a singular point of the Legendre curve $(\gamma,\nu)$, that is, $(\ell(t_0),\beta(t_0))=(0,0)$. 

\begin{definition}\label{vertex.def}{\rm
We say that $t_0$ is a {\it vertex of the frontal $\gamma$} (or, {\it of the Legendre curve $(\gamma,\nu)$}) if the evolute $\mathcal{E}v(\gamma)$ of the frontal exists and $\dot{\mathcal{E}}v(\gamma)(t_0)=0$.
}
\end{definition}

Suppose that the evolute of the frontal $\mathcal{E}v(\gamma)$ exists and $\beta=\alpha \ell$. 
Then $t_0$ is a vertex of the frontal $\gamma$ if and only if $\dot{\alpha}(t_0)=0$.

\begin{remark}{\rm 
Suppose that $t_0$ is a singular point $(\gamma,\nu)$. 
If $t_0$ is a vertex of $\gamma$, then $\gamma$ and the osculating circle have at least fourth order contact at $t_0$, up to congruence as Legendre curves (cf. \cite{Fukunaga-Takahashi-2015}). 
}
\end{remark}

\begin{proposition}[\cite{Fukunaga-Takahashi-2015}]\label{evolute-parallel.frontal}
Let $(\gamma, \nu) : I \to \mathbb{R}^2 \times S^1$ be a Legendre curve with curvature $(\ell, \beta)$. 
If the evolute $\mathcal{E}v(\gamma)$ exists, then the evolute of a parallel curve of $\gamma$ exists. 
Moreover, the evolute $\mathcal{E}v(\gamma^\lambda)$ coincides with the evolute  $\mathcal{E}v(\gamma)$. 
\end{proposition}

\begin{proposition}
Let $(\gamma, \nu) : I \to \mathbb{R}^2 \times S^1$ be a Legendre curve with the curvature $(\ell, \beta)$ and $\beta=\alpha \ell$. 
Suppose that $t:\widetilde{I} \to I$ is a parameter change and $(\widetilde{\gamma},\widetilde{\nu})=(\gamma \circ t,\nu \circ t)$, that is, $(\widetilde{\gamma},\widetilde{\nu})$ and $(\gamma,\nu)$ are $\mathcal{R}$-equivalent. 
Then $u_0 \in \widetilde{I}$ is a vertex of $\widetilde{\gamma}$ if and only if $t(u_0) \in I$ is a vertex of $\gamma$.
\end{proposition}
\demo
By Proposition \ref{Legendre.function}, the curvature of $(\widetilde{\gamma},\widetilde{\nu})$ is given by $(\widetilde{\ell},\widetilde{\beta})=((\ell \circ t) t', (\beta \circ t) t')$.
It follows that $\widetilde{\alpha}=\alpha \circ t$ and $\mathcal{E}v(\widetilde{\gamma})=\mathcal{E}v(\gamma) \circ t$, that is,  $\widetilde{\alpha}$ (respectively, $\mathcal{E}v(\widetilde{\gamma})$) and $\alpha$ (respectively, $\mathcal{E}v(\gamma)$) are $\mathcal{R}$-equivalent. 
Thus, $u_0 \in \widetilde{I}$ is a vertex of $\widetilde{\gamma}$ if and only if $t(u_0) \in I$ is a vertex of $\gamma$.
\enD

\begin{proposition}\label{affine1}
Let $(\gamma,\nu):I \to \R^2 \times S^1$ be a Legendre curve with curvature $(\ell,\beta)$ and $\Phi:\R^2 \to \R^2$ be a diffeomorphism. 
We denote $\nu=(a,b)$ and $\widetilde{\gamma}=\Phi \circ \gamma$.
Then we have the following.
\par
$(1)$ Suppose that the diffeomorphism $\Phi:\R^2 \to \R^2$ is given by  $\Phi(x,y)=(a_{11}x+a_{12}y,a_{21}x+a_{22}y)$, where $a_{11}a_{22}-a_{12}a_{21} \not=0$ and $a_{11},a_{12},a_{21},a_{22} \in \R$. 
Then $t_0$ is a singular point and a vertex of $\gamma$ if and only if $t_0$ is a singular point and a vertex of $\widetilde{\gamma}=\Phi \circ \gamma$.
\par
$(2)$ Suppose that the diffeomorphism $\Phi: \R^2 \to \R^2$ is given by $\Phi(x,y)=(y,x)$. 
Then $t_0$ is a vertex of $\gamma$ if and only if $t_0$ is a vertex of $\widetilde{\gamma}=\Phi \circ \gamma$.
\end{proposition}
\demo 
$(1)$ By Corollary \ref{change-target} $(1)$, the curvature of $(\widetilde{\gamma},\widetilde{\nu}):I \to \R^2 \times S^1$ is given by $(\widetilde{\ell},\widetilde{\beta})=\left((a_{11}a_{22}-a_{12}a_{21})\ell/|\overline{\nu}|^2,|\overline{\nu}|\beta\right)$, where $\widetilde{\nu}=\overline{\nu}/|\overline{\nu}|$ and $\overline{\nu}=(a_{22}a-a_{21}b,-a_{12}a+a_{11}b)$. 
Therefore, $t_0$ is also a singular point of $\widetilde{\gamma}$.  
Since $\widetilde{\beta}=(|\overline{\nu}|^3 \alpha/(a_{11}a_{22}-a_{12}a_{21}))\widetilde{\ell}$, 
we have $\widetilde{\alpha}=|\overline{\nu}|^3 \alpha/(a_{11}a_{22}-a_{12}a_{21})$. 
By $\beta(t_0)=0$, we have $\alpha(t_0)=0$ or $\ell(t_0)=0$.
By a direct calculation, 
\begin{align*}
\dot{\widetilde{\alpha}}(t)&=\frac{1}{a_{11}a_{22}-a_{12}a_{21}}\Bigl(3|\overline{\nu}(t)|((-a_{22}b(t)-a_{21}a(t))(a_{22}a(t)-a_{21}b(t))\\
&\quad +(a_{12}b(t)+a_{11}a(t))(-a_{12}a(t)+a_{11}b(t)))\alpha(t)\ell(t)+|\overline{\nu}(t)|^3\dot{\alpha}(t) \Bigr).
\end{align*}
If $\dot{\alpha}(t_0)=0$, then $\dot{\widetilde{\alpha}}(t_0)=0$. 
It follows that $t_0$ is also a vertex of $\widetilde{\gamma}$. 
The converse is follows from $\Phi^{-1}$.
\par
$(2)$ By Corollary \ref{change-target} $(2)$, the curvature of $(\widetilde{\gamma},\widetilde{\nu}):I \to \R^2 \times S^1$ is given by $(\widetilde{\ell},\widetilde{\beta})=\left(-\ell,\beta\right)$, where $\widetilde{\nu}=(-b,-a)$. 
Since $\widetilde{\beta}=-\alpha \widetilde{\ell}$, we have $\widetilde{\alpha}=-\alpha$. 
Therefore, $t_0$ is a vertex of $\gamma$ if and only if $t_0$ is a vertex of $\widetilde{\gamma}$.
\enD

We investigate the four vertex theorem of frontals. 
We give conditions that $\gamma$ has at least four vertices by using the curvature of the Legendre curve and type of singular points of the Legendre curve.
By Roll's theorem, we have the following result.
\begin{proposition}\label{condition1}
Let $(\gamma, \nu) : I=[a,b] \to \mathbb{R}^2 \times S^1$ be a closed Legendre curve with curvature $(\ell, \beta)$ and $\beta=\alpha \ell$.
\par
$(1)$ If there exist at least two points such that $\alpha=\dot{\alpha}=0$ at the points, then $\gamma$ has at least four vertices. 
\par
$(2)$ If $\alpha$ has at least four zero points, then $\gamma$ has at least four vertices. 
\end{proposition}

Suppose that a Legendre curve $(\gamma,\nu):I \to \R^2 \times S^1$ at $t_0 \in I$ is of type $(n,m)$, that is, 
$(\gamma,\nu)$ at $t_0$ is $\mathcal{R}$-equivalent to $(\widetilde{\gamma},\widetilde{\nu}):(\R,0) \to \R^2 \times S^1$,
\begin{align*}
\widetilde{\gamma}(t) =(\pm t^n,t^mf(t)),\
\widetilde{\nu}(t) =\frac{(-mt^kf(t)-t^{k+1}\dot{f}(t), \pm n)}{\sqrt{(mt^kf(t)+t^{k+1}\dot{f}(t))^2+n^2}},
\end{align*}
where $f:(\R,0) \to \R$ is a smooth function germ with $f(0) \not=0$ and the curvature $(\widetilde{\ell},\widetilde{\beta})$ (see Example \ref{nm-type}). 
Then we say that {\it $\gamma$ has a point $t_0$ of type $(n,m)$}.
\par
If $n \ge k$, then there exists a unique function $\widetilde{\alpha}:(\R,0) \to \R$, 
$$
\widetilde{\alpha}(t)=\mp \frac{t^{n-k}((mt^kf(t)+t^{k+1}\dot{f}(t))^2+n^2)^{\frac{3}{2}}}{n(mkf(t)+(m+k+1)t\dot{f}(t)+t^{2}\ddot{f}(t))}
$$
such that $\widetilde{\beta}=\widetilde{\alpha} \widetilde{\ell}$.

\begin{proposition}\label{condition2}
Let $(\gamma, \nu) : I \to \mathbb{R}^2 \times S^1$ be a Legendre curve with curvature $(\ell, \beta)$ and $\beta=\alpha \ell$. 
\par
$(1)$ In the case of $n=k$. 
If $\gamma$ has four points of type $(n,m)=(n,2n)$ with $\dot{f}(0)=0$, then $\gamma$ has at least four vertices. 
\par
$(2)$ In the case of $n=k+1$. 
If $\gamma$ has four points of type $(n,m)=(n,2n-1)$, then $\gamma$ has at least four vertices. 
\par
$(3)$ In the case of $n \ge k+2$. 
If $\gamma$ has two points of type $(n,m)$, then $\gamma$ has at least four vertices. 
\end{proposition}
\demo
If $\gamma$ has a point of type $(n,m)$, then $\alpha$ at the point is  $\mathcal{R}$-equivalent to $\widetilde{\alpha}$ at $0$. 
\par
$(1)$ If $n=k$, $\widetilde{\alpha}$ is given by
$$ 
\widetilde\alpha(t)=\mp \frac{((mt^nf(t)+t^{n+1}\dot{f}(t))^2+n^2)^{\frac{3}{2}}}{n(mnf(t)+(m+n+1)t\dot{f}(t)+t^{2}\ddot{f}(t))}.
$$
By a direct calculation, we have 
\begin{align*}
\dot{\widetilde{\alpha}}(t)&=\mp \frac{\sqrt{(mt^nf(t)+t^{n+1}\dot{f}(t))^2+n^2}}{n(mnf(t)+(m+n+1)t\dot{f}(t)+t^{2}\ddot{f}(t))^2}\\
&\Bigl(((mt^nf(t)+t^{n+1}\dot{f}(t))^2+n^2)((mn+m+n+1)\dot{f}(t)+(m+n+3)t\ddot{f}(t)+t^2\dddot{f}(t))\\
&\quad -3t^{n}(mf(t)+t\dot{f}(t))(mnf(t)+(m+n+1)t\dot{f}(t)+t^2\ddot{f}(t))^2\Bigr).
\end{align*}
By $\dot{f}(0)=0$, we have $\dot{\widetilde{\alpha}}(0)=0$. 
Since $\alpha$ and $\widetilde{\alpha}$ are $\mathcal{R}$-equivalent, if there exist   four points $t_1, t_2, t_3$ and $t_4 \in I$ of $\gamma$ of type $(n,2n)$ with $\dot{f}(0)=0$, then $\dot{\alpha}(t_1)=\dot{\alpha}(t_2)=\dot{\alpha}(t_3)=\dot{\alpha}(t_4)=0$. 
Therefore, $\gamma$ has at least four vertices. 
\par
$(2)$ If $n=k+1$, $\widetilde{\alpha}$ is given by
$$
\widetilde\alpha(t)=\mp \frac{t((mt^{n-1}f(t)+t^{n}\dot{f}(t))^2+n^2)^{\frac{3}{2}}}{n(m(n-1)f(t)+(m+n)t\dot{f}(t)+t^{2}\ddot{f}(t))},
$$
$\widetilde{\alpha}(0)=0$. 
If there exist four points $t_1, t_2, t_3$ and $t_4\in I$ of $\gamma$ of type $(n,2n-1)$, $\alpha(t_1)=\alpha(t_2)=\alpha(t_3)=\alpha(t_4)=0$. 
By Proposition \ref{condition1} (2), $\gamma$ has at least four vertices. 
\par
$(3)$ If $n \ge k+2$, $\widetilde{\alpha}$ and $\dot{\widetilde{\alpha}}$ are given by
\begin{align*}
\widetilde\alpha(t)&=\mp \frac{t^{n-k}((mt^kf(t)+t^{k+1}\dot{f}(t))^2+n^2)^{\frac{3}{2}}}{n(mkf(t)+(m+k+1)t\dot{f}(t)+t^{2}\ddot{f}(t))}, \\
\dot{\widetilde\alpha}(t)&=\mp \frac{t^{n-k-1}\sqrt{(mt^kf(t)+t^{k+1}\dot{f}(t))^2+n^2}}{n(mkf(t)+(m+k+1)t\dot{f}(t)+t^{2}\ddot{f}(t))^2}\\
&\Bigl(\bigl((mt^kf(t)+t^{k+1}\dot{f}(t))^2+n^2\bigr)\bigl((n-k)(mkf(t)+(m+k+1)t\dot{f}(t)+t^2\ddot{f}(t))\\
&\hspace{5mm}-(mk+m+k+1)t\dot{f}(t)-(m+k+3)t^2\ddot{f}(t)+t^3\dddot{f}(t)\bigr)\\
&\hspace{10mm}-(mf(t)+t\dot{f}(t))(mkf(t)+(m+k+1)t\dot{f}(t)+t^2\ddot{f}(t))^2\Bigr),
\end{align*}
$\widetilde{\alpha}(0)=\dot{\widetilde{\alpha}}(0)=0$. 
If there exist two points $t_1$ and $t_2\in I$ of $\gamma$ of type $(n,m)$, $\alpha(t_1)=\dot\alpha(t_1)=0$ and $\alpha(t_2)=\dot\alpha(t_2)=0$. 
By Proposition \ref{condition1} (1), $\gamma$ has at least four vertices. 
\enD
\begin{theorem}\label{four-vertex-theorem-frontal}
Let $(\gamma, \nu) : I=[a,b] \to \mathbb{R}^2 \times S^1$ be a closed Legendre curve with curvature $(\ell, \beta)$ and $\beta=\alpha \ell$. 
Suppose that zeros of $\ell$ and of $\beta$ are isolated points.
If $\gamma$ is a simple convex frontal, then $\gamma$ has at least four vertices.
\end{theorem}
\demo
By Theorem \ref{Convex-frontal}, we may consider the case $(i)$ $\ell(t) \ge 0$ and $\beta(t) \ge 0$ for all $t \in I$. 
We can similarly prove the other cases $(ii)$, $(iii)$ and $(iv)$.
\par
Since $\beta(t)=\alpha(t) \ell(t)$, $\alpha(t) \ge 0$ for all $t \in I$.
By Proposition \ref{evolute-parallel.frontal}, the evolute of parallel curves of the Legendre curve $\mathcal{E}v(\gamma^\lambda)$ is the same as the evolute of the frontal $\mathcal{E}v(\gamma)$. 
Hence vertices of the parallel curves are the same as the vertices of $\gamma$.
We consider a parallel curve $(\gamma^\lambda,\nu)$ for a constant $\lambda>0$. 
By Proposition \ref{parallel.frontal}, the curvature $(\ell^\lambda,\beta^\lambda)$ of $(\gamma^\lambda,\nu)$ is given by $(\ell,\beta+\lambda \ell)$.
It follows that $\alpha^\lambda=\alpha+\lambda$ is a positive function. 
Since $\gamma$ is a simple convex frontal and the rotation index of Legendre curve $(\gamma,\nu)$ is $\pm 1$ by \cite[Lemma 2.5]{Fukunaga-Takahashi-2016}, we can show that rotation index of Legendre curve $(\gamma^\lambda,\nu)$ is the same as the rotation index of Legendre curve $(\gamma,\nu)$.
Moreover, the signs of $\ell, \beta$ and $\beta^\lambda$ do not change, $\gamma^\lambda$ is also a convex frontal (cf. \cite[Lemma 2.6]{Fukunaga-Takahashi-2016}). 
Since $(\gamma,\nu)$ is a closed Legendre curve, $(\gamma^\lambda,\nu)$ is also a closed Legendre curve. 
We denote $\gamma^\lambda(t)=(x^\lambda(t), y^\lambda(t))$ and $\nu(t)=(a(t),b(t))$.  
Then we have $\mu(t)=(-b(t),a(t))$. 
\par
\begin{figure}[htbp]
  \centering
      \includegraphics[scale=0.17]{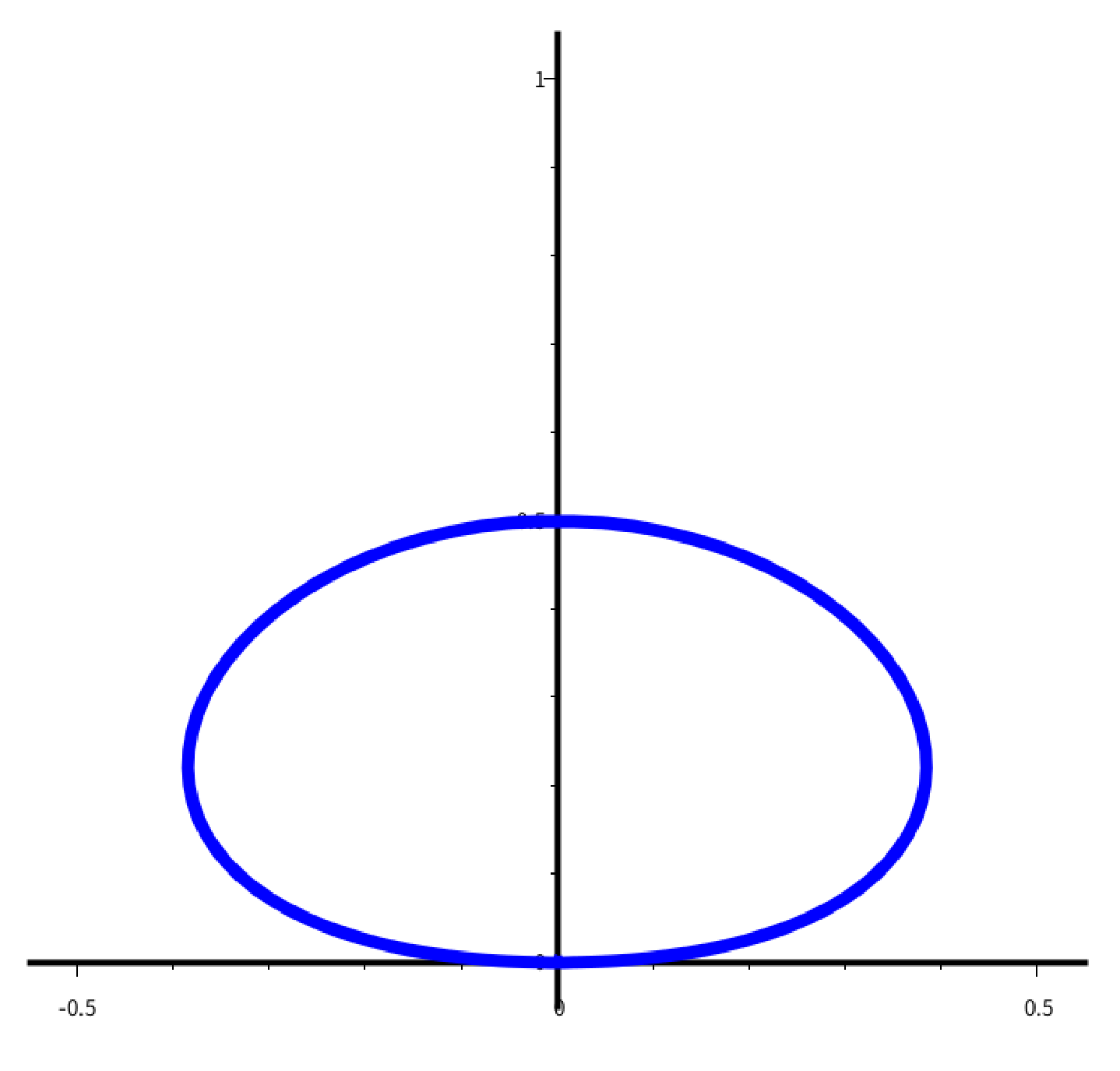}
      \includegraphics[scale=0.17]{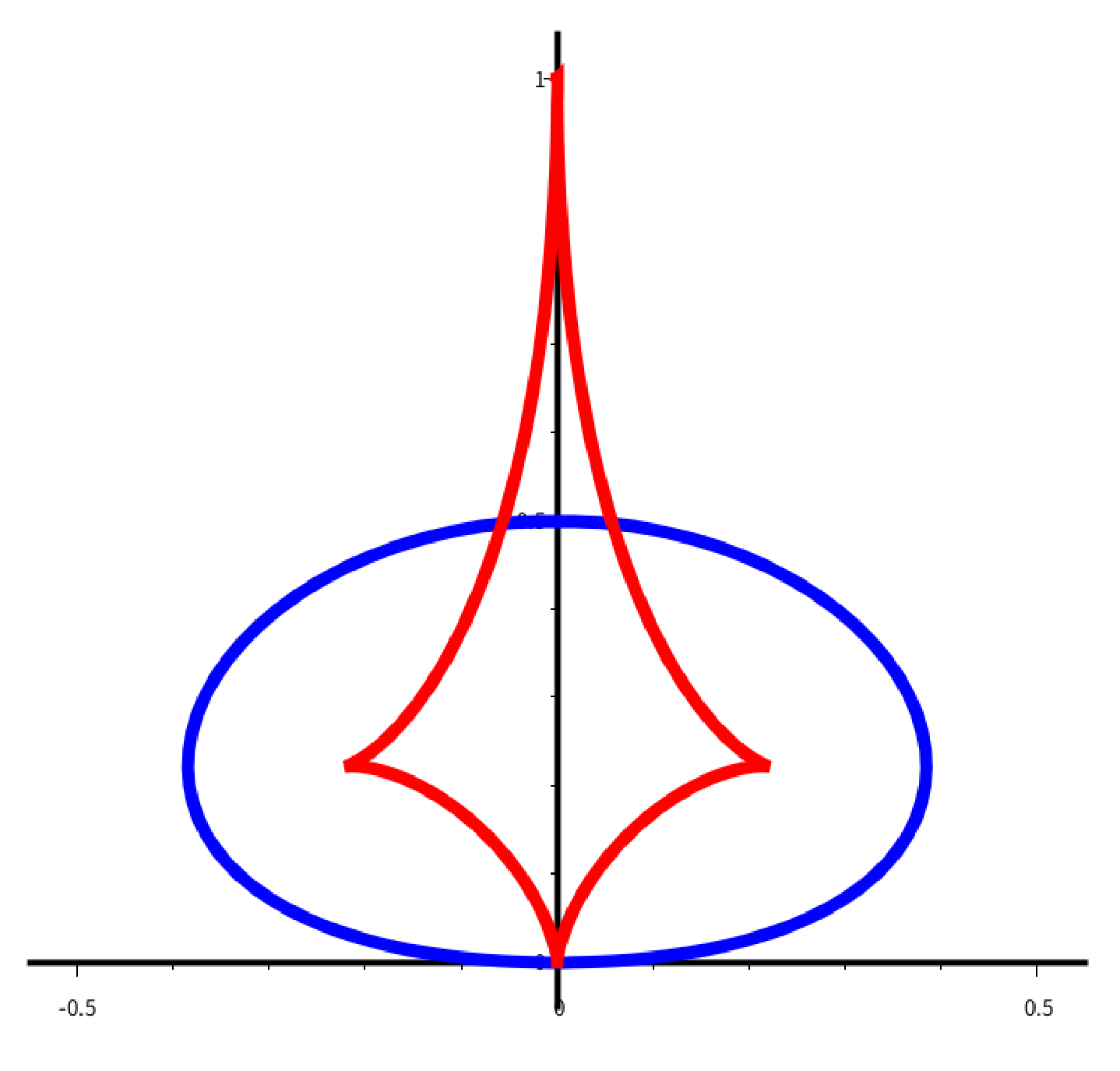}
  \caption{\centering{$\gamma(t)$ and $\mathcal{E}v(\gamma)(t)$}}
  \label{figure1}
\end{figure}
\begin{figure}[htbp]
  \centering
      \includegraphics[scale=0.17]{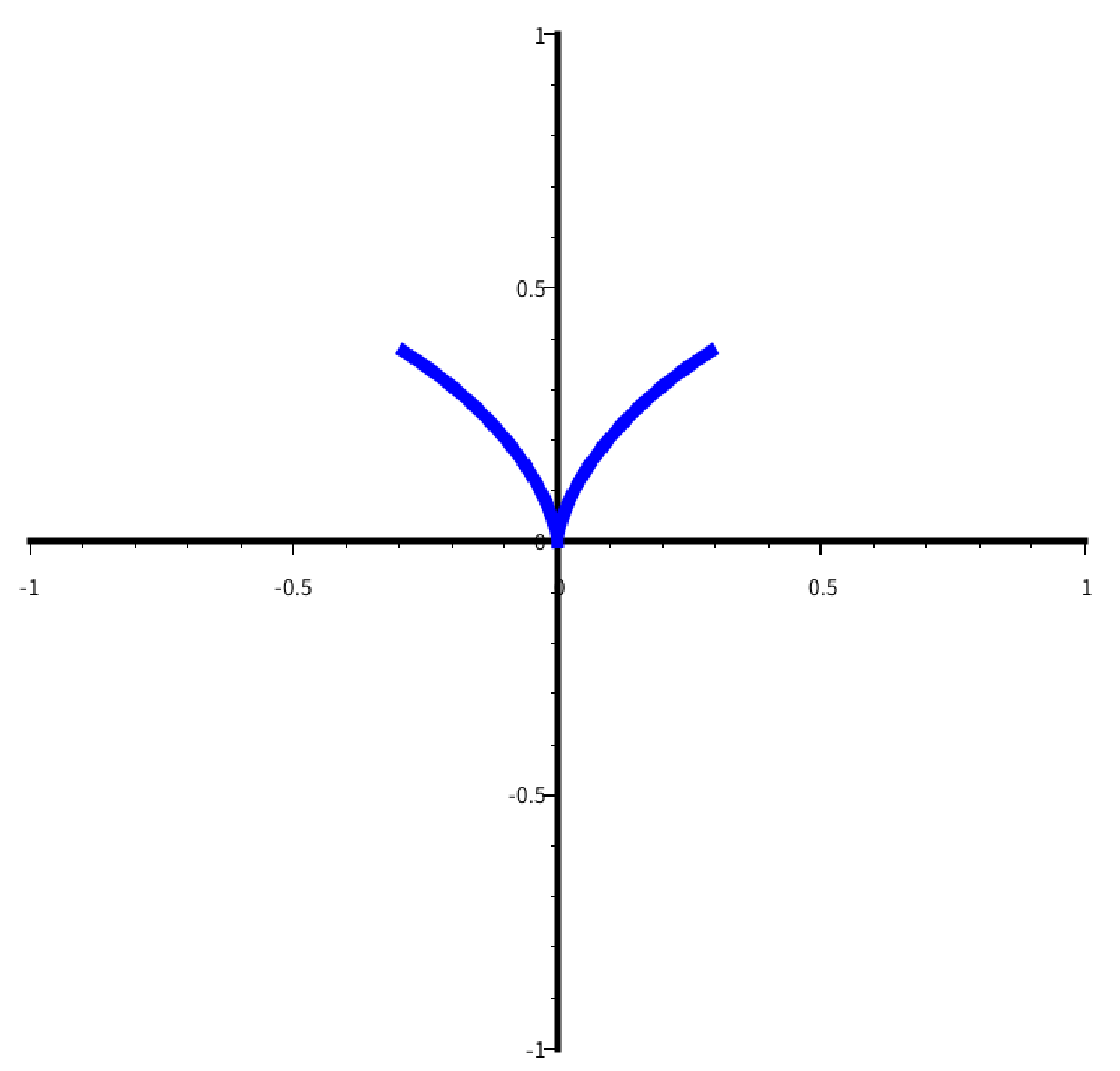}
      \includegraphics[scale=0.17]{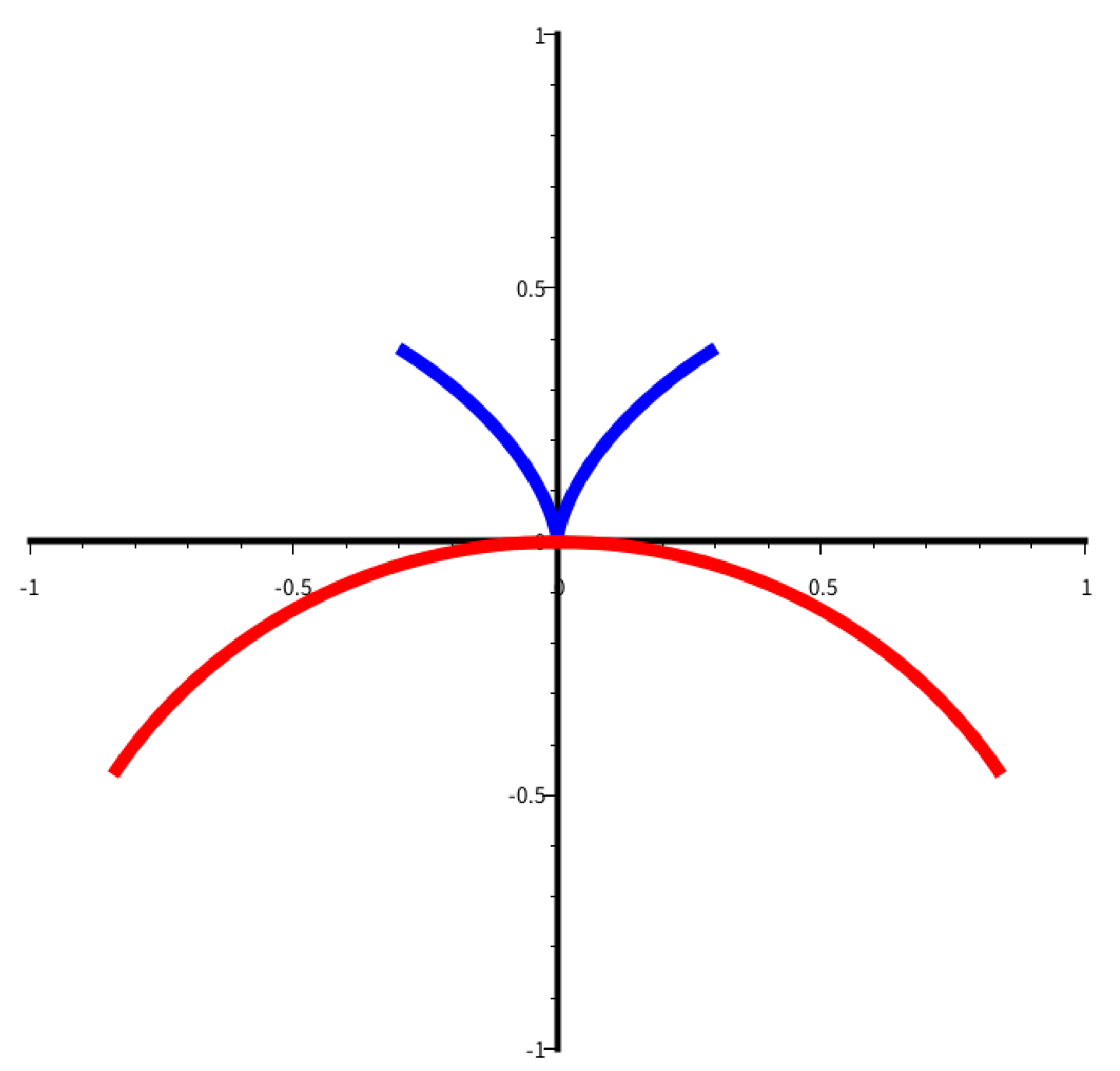}
  \caption{\centering{$\gamma(t)$ and $\mathcal{E}v(\gamma)(t)$}}
  \label{figure2}
\end{figure}

Now suppose that $\gamma$ (and also $\gamma^\lambda$) has two vertices. 
By using congruence, we may 
\begin{align*}
\int^b_a \frac{d}{dt}\left(\frac{1}{\alpha^\lambda(t)}\right) y^{\lambda}(t) dt= -\int^b_a \frac{\dot{\alpha}^\lambda(t)}{\alpha^{\lambda}(t)^2} y^{\lambda}(t) dt\not=0,
\end{align*}
see the case of regular curves in \cite{Gibson}. 
On the other hand, 
\begin{align*}
&\int^b_a \frac{d}{dt} \left(\frac{1}{\alpha^\lambda(t)}\right) y^{\lambda}(t) dt
=\left[ \left(\frac{1}{\alpha^\lambda(t)}\right) y^{\lambda}(t) \right]^b_a
-\int^b_a \left(\frac{1}{\alpha^\lambda(t)}\right) \dot{y}^{\lambda}(t) dt \\
=&-\int^b_a \left(\frac{1}{\alpha^\lambda(t)}\right) \beta^\lambda(t) a(t) dt
=-\int^b_a \ell(t) a(t) dt=-[b(t)]^b_a=0.
\end{align*}
This is a contradiction. 
It follows that $\gamma$ has at least three vertices. 
By the same arguments, we can show that $\gamma$ has at least four vertices. 
\enD
\begin{example}
{\rm
Let $(\gamma,\nu):[0,\pi] \to \R^2 \times S^1$ be 
\begin{align*}
 \gamma(t)&=\left(\sin^2 t \cos t,\frac{1}{2}\sin^4 t \right)=\left(\frac{1}{2}\sin 2t\sin t,\frac{1}{2}\sin^2 t-\frac{1}{8}\sin^2 2t\right) ,\\
 \nu(t)&=\frac{1}{\sqrt{(\cos 2t+\cos^2 t)^2+\cos^2 t(1-\cos 2t)^2}}\left(-\cos t(1-\cos 2t), \cos 2t+\cos^2 t\right).
\end{align*}
By a direct calculation, $(\gamma,\nu)$ is a closed Legendre curve with the curvature
\begin{align*}
\ell(t)&=\frac{2\sin t \left((\cos 2t+\cos^2 t)^2+3\cos^2 t(1-\cos 2t)\right)}{(\cos 2t+\cos^2 t)^2+\cos^2 t(1-\cos 2t)^2 }, \\ 
\beta(t)&=-\sin t\sqrt{(\cos 2t+\cos^2 t)^2+\cos^2 t(1-\cos 2t)^2}.
\end{align*}
Since $\ell$ is non-negative and $\beta$ is non-positive on $[0,\pi]$, $\gamma$ is a simple convex closed frontal. 
There exists a unique function $\alpha:[0,\pi] \rightarrow \R$, 
\begin{align*}
\alpha(t)=-\frac{\left((\cos 2t+\cos^2 t)^2+\cos^2 t(1-\cos 2t)^2 \right)^\frac{3}{2}}{2\left((\cos 2t+\cos^2 t)^2+3\cos^2 t(1-\cos 2t)\right)}
\end{align*}
such that $\beta=\alpha\ell$.
Then the evolute of the frontal is given by
\begin{align*}
&\mathcal{E}v(\gamma)(t)=\Biggl(\frac{1}{2}\sin 2t\sin t-\frac{\cos t(1-\cos 2t)\left((\cos 2t+\cos^2 t)^2+\cos^2 t(1-\cos 2t)^2\right)}{2\left((\cos 2t+\cos^2 t)^2+3\cos^2 t(1-\cos 2t)\right)}, \\
&\qquad \frac{1}{2}\sin^2 t-\frac{1}{8}\sin^2 2t+\frac{(\cos 2t+\cos^2 t)\left((\cos 2t+\cos^2 t)^2+\cos^2 t(1-\cos 2t)^2\right)}{2\left((\cos 2t+\cos^2 t)^2+3\cos^2 t(1-\cos 2t)\right)}\Biggr).
\end{align*}
Since 
\begin{align*}
\dot{\alpha}(t)&=-\frac{3\sin 2t\left((\cos 2t+\cos^2 t)^2+\cos^2 t(1-\cos 2t)^2\right)^\frac{1}{2}}{2\left((\cos 2t+\cos^2 t)^2+3\cos^2 t(1-\cos 2t)\right)^2}g(t),
\end{align*} 
where $g(t)=-10\cos^8 t-\cos^6 t-7\cos^4 t+\cos^2 t+1$,
$\dot{\alpha}(0)=\dot{\alpha}(\pi/2)=0$. It follows that $t=0$ and $\pi/2$ are vertices of $\gamma$. 
Moreover, there exist only two points $t_1\in[0, \pi/2)$ and $t_2\in[\pi/2, \pi)$ such that $g(t_1)=g(t_2)=0$. 
Therefore, $\gamma$ has four vertices. See Figure. \ref{figure1}.
}
\end{example}
\begin{example}{\rm
Let $(\gamma,\nu):[0,2\pi] \to \R^2 \times S^1$ be 
\begin{align*}
 \gamma(t)&=(\sin t\cos(\sin t)-\sin(\sin t), \sin t\sin(\sin t)+\cos(\sin t)-1) ,\\
 \nu(t)&=(\cos(\sin t), \sin(\sin t)).
\end{align*}
By a direct calculation, $(\gamma,\nu)$ is a closed Legendre curve with the curvature
$
\ell(t)=\cos t, \ \beta(t)=\cos t\sin t.
$
It follows that there exists a unique function $\alpha:[0,2\pi] \rightarrow \R$, $\alpha(t)=\sin t$ such that $\beta=\alpha\ell$. 
Then the evolute of the frontal $\gamma$ is given by 
\begin{align*}
\mathcal{E}v(\gamma)(t)&=\gamma(t)-\alpha(t)\nu(t)=(-\sin(\sin t), \cos(\sin t)-1).
\end{align*}
Since $\dot{\alpha}(t)=\cos t$, $\gamma$ has only two vertices. 
Note that $\gamma$ is neither a simple curve nor a convex curve.
Therefore, it does not satisfy on the four vertex theorem. See Figure. \ref{figure2}.
}
\end{example}

\section{Vertices of frontals in the general case}

In the previous section, we define the vertex using the evolute of the Legendre curve $(\gamma,\nu):I \to \R^2 \times S^1$ with curvature $(\ell,\beta)$ under the condition $\beta=\alpha \ell$. 
It follows that $\ell$ and $\beta$ are linearly dependent on $I$. 
Therefore, we give a definition of a generalisation of the vertex of frontals.

Let $(\gamma,\nu):I \to \R^2 \times S^1$ be a Legendre curve with curvature $(\ell,\beta)$.
In this section, we assume that $\ell$ and $\beta$ are linearly dependent on $I$, that is, there exists a non-zero smooth mapping $(k_1,k_2):I \to \R^2 \setminus \{0\}$ such that $k_1(t) \ell(t)+k_2(t) \beta(t)=0$ for all $t \in I$.
\begin{definition}\label{generalised-vertex}{\rm 
Under the above assumption, we say that $t_0 \in I$ is a {\it vertex} of the frontal $\gamma$ (or, {\it of the Legendre curve $(\gamma,\nu)$}) if 
$\dot{k}_1(t_0) k_2(t_0)-k_1(t_0) \dot{k}_{2}(t_0)=0$.
}
\end{definition}

If $\beta(t) \not=0$ for all $t \in I$, then $\gamma$ is a regular curve and we can take $k_1(t)=1, k_2(t)=-\ell(t)/\beta(t)$. 
It follows that $t_0$ is a vertex of $\gamma$ as a regular curve if and only if  $t_0 \in I$ is a vertex of $\gamma$ as the frontal. 
Moreover, if $\ell(t) \not=0$ for all $t \in I$, then $(\gamma,\nu)$ is a Legendre immersion and we can take $k_1(t)=-\beta(t)/\ell(t), k_2(t)=1$. 
It follows that $t_0$ is a vertex of the front $\gamma$ if and only if  $t_0 \in I$ is a vertex of $\gamma$ as the frontal. 
Furthermore, if $\beta(t)=\alpha(t) \ell(t)$ for all $t \in I$, we can take $k_1(t)=-\alpha(t), k_2(t)=1$. 
It follows that $t_0$ is a vertex of the frontal $\gamma$ in the sense of the previous section if and only if  $t_0 \in I$ is a vertex of $\gamma$ as the frontal. 


\begin{proposition}\label{evolute-parallel2.frontal}
If $t_0 \in I$ is a vertex of $\gamma$, then $t_0$ is also a vertex of the parallel curve $\gamma^\lambda$.
\end{proposition}
\demo 
By Proposition \ref{parallel.frontal}, the curvature $(\ell^\lambda, \beta^\lambda)$ of $(\gamma^\lambda, \nu)$ is given by $(\ell^\lambda, \beta^\lambda)=(\ell, \beta+\lambda\ell)$. 
We consider a smooth mapping $(k_1^\lambda, k_2^\lambda):I \to \R^2$ such that
$$
\left(
\begin{array}{c}
k_1^\lambda(t)\\
k_2^\lambda(t)
\end{array}
\right)
=
\left(
\begin{array}{cc}
1 & -\lambda\\
0 & 1
\end{array}
\right)
\left(
\begin{array}{c}
k_1(t)\\
k_2(t)
\end{array}
\right). 
$$
Since $(k_1,k_2):I \to \R^2 \setminus \{0\}$ is a non-zero smooth mapping such that $k_1(t) \ell(t)+k_2(t) \beta(t)=0$ for all $t \in I$, $(k_1^\lambda, k_2^\lambda)$ is also a non-zero smooth mapping and 
\begin{align*}
k_1^\lambda(t)\ell^\lambda(t)+k_2^\lambda(t)\beta^\lambda(t)&=(k_1(t)-\lambda k_2(t))\ell(t)+k_2(t)(\beta(t)+\lambda\ell(t))\\\
&=k_1(t)\ell(t)+k_2(t)\beta(t)=0.
\end{align*}
By a direct calculation, 
\begin{align*}
\dot{k_1^\lambda}(t) k_2^\lambda(t)-k_1^\lambda(t)\dot{k_2^\lambda}(t)&=(\dot{k}_1(t)-\lambda\dot{k}_2(t))k_2(t)-(k_1(t)-\lambda k_2(t))\dot{k}_2(t)\\
&=\dot{k}_1(t) k_2(t)-k_1(t) \dot{k}_{2}(t).
\end{align*}
Therefore, if $t_0 \in I$ is a vertex of $\gamma$, then $t_0$ is also a vertex of $\gamma^\lambda$.
\enD
By the symmetry of the linearly dependent condition $k_1(t) \ell(t)+k_2(t) \beta(t)=0$ for $(\ell,\beta)$, we have the following.
\begin{proposition}\label{evolute-parallel3.frontal}
Let $(\widetilde{\gamma}^\lambda,\widetilde{\nu}^\lambda):I \to \R^2 \times S^1$ be a Legendre curve with the curvature $(\widetilde{\ell}^\lambda,\widetilde{\beta}^\lambda)=(\ell+\lambda \beta,\beta)$ for $\lambda \in \R$.
If $t_0 \in I$ is a vertex of $\gamma$, then $t_0$ is also a vertex of $\widetilde{\gamma}^\lambda$.
\end{proposition}
\demo
We consider a smooth mapping $(\widetilde{k}_1^\lambda, \widetilde{k}_2^\lambda):I \to \R^2$ such that 
$$
\left(
\begin{array}{c}
\widetilde{k}_1^\lambda(t)\\
\widetilde{k}_2^\lambda(t)
\end{array}
\right)
=
\left(
\begin{array}{cc}
1 & 0\\
-\lambda & 1
\end{array}
\right)
\left(
\begin{array}{c}
k_1(t)\\
k_2(t)
\end{array}
\right).
$$
Since $(k_1,k_2):I \to \R^2 \setminus \{0\}$ is a non-zero smooth mapping such that $k_1(t) \ell(t)+k_2(t) \beta(t)=0$ for all $t \in I$, $(\widetilde{k}_1^\lambda, \widetilde{k}_2^\lambda)$ is also a non-zero smooth mapping and 
\begin{align*}
\widetilde{k}_1^\lambda(t)\widetilde{\ell}^\lambda(t)
+\widetilde{k}_2^\lambda(t)\widetilde{\beta}^\lambda(t)
&=k_1(t)(\ell(t)+\lambda\beta(t))+(k_2(t)-\lambda k_1(t))\beta(t)\\
&=k_1(t)\ell(t)+k_2(t)\beta(t)=0.
\end{align*}
By a direct calculation, 
\begin{align*}
\dot{\widetilde{k}_1^\lambda}(t) \widetilde{k}_2^\lambda(t)-\widetilde{k}_1^\lambda(t)\dot{\widetilde{k}_2^\lambda}(t)&=\dot{k}_1(t)(k_2(t)-\lambda k_1(t))-k_1(t)(\dot{k}_2(t)-\lambda \dot {k}_1(t))\\
&=\dot{k}_1(t) k_2(t)-k_1(t) \dot{k}_{2}(t).
\end{align*}
Therefore, if $t_0 \in I$ is a vertex of $\gamma$, then $t_0$ is also a vertex of $\widetilde{\gamma}^\lambda$.
\enD

More generally, we have the following.

\begin{proposition}\label{evolute-affine.frontal}
Let $(\gamma[A],\nu[A]):I \to \R^2 \times S^1$ be a Legendre curve with the curvature $({\ell}[A],{\beta}[A])=A ^t(\ell,\beta)=(a_{11} \ell+a_{12} \beta, a_{21} \ell+a_{22} \beta)$, where
$$
A=\begin{pmatrix}
a_{11} & a_{12} \\
a_{21} & a_{22}
\end{pmatrix} \in GL(2,\R).
$$
Then $t_0 \in I$ is a vertex of $\gamma$ if and only if $t_0$ is a vertex of ${\gamma}[A]$.
\end{proposition}
\demo 
We consider a smooth mapping $({k_1[A]}, {k_2[A]}):I \to \R^2$ such that 
\begin{align*}
\left(
\begin{array}{c}
k_1[A](t)\\
k_2[A](t)
\end{array}
\right)
&=\det A {(^{t}A)}^{-1}
\left(
\begin{array}{c}
k_1(t)\\
k_2(t)
\end{array}
\right)
=
\left(
\begin{array}{cc}
a_{22} & -a_{21}\\
-a_{12} & a_{11}
\end{array}
\right)
\left(
\begin{array}{c}
k_1(t)\\
k_2(t)
\end{array}
\right).
\end{align*}
Since $(k_1,k_2):I \to \R^2 \setminus \{0\}$ is a non-zero smooth mapping such that $k_1(t) \ell(t)+k_2(t) \beta(t)=0$ for all $t \in I$, $({k_1[A]}, {k_2[A]})$ is also a non-zero smooth mapping and 
\begin{align*}
{k_1[A]}(t){\ell[A]}(t)+{k_2[A]}(t){\beta[A]}(t)&=(a_{11}a_{22}-a_{12}a_{21})(k_1(t)\ell(t)+k_2(t)\beta(t))=0.
\end{align*}
By a direct calculation, 
\begin{align*}
\dot{k_1}[A](t) {k_2[A]}(t)-{k_1[A]}(t)\dot{{k_2}}[A](t)&=(a_{11}a_{22}-a_{12}a_{21})(\dot{k}_1(t) k_2(t)-k_1(t) \dot{k}_{2}(t)).
\end{align*}
Therefore,  $t_0 \in I$ is a vertex of $\gamma$ if and only if $t_0$ is a vertex of ${\gamma}[A]$.
\enD

\begin{proposition}
Suppose that $t:\widetilde{I} \to I$ is a parameter change and $(\widetilde{\gamma},\widetilde{\nu})=(\gamma \circ t,\nu \circ t)$, that is, $(\widetilde{\gamma},\widetilde{\nu})$ and $(\gamma,\nu)$ are $\mathcal{R}$-equivalent. 
Then $u_0 \in \widetilde{I}$ is a vertex of $\widetilde{\gamma}$ if and only if $t(u_0) \in I$ is a vertex of $\gamma$.
\end{proposition}
\demo
By Proposition \ref{Legendre.function}, the curvature $(\widetilde\ell, \widetilde\beta)$ of $(\widetilde\gamma, \widetilde\nu)$ is given by $(\widetilde\ell, \widetilde\beta)=((\ell \circ t)t', (\beta \circ t)t')$. 
We consider a smooth mapping $(\widetilde{k}_1, \widetilde{k}_2): \widetilde{I} \to \R^2$, $(\widetilde{k}_1, \widetilde{k}_2)(u)=(k_1\circ t(u), k_2 \circ t(u))$. 
Since $(k_1,k_2):I \to \R^2 \setminus \{0\}$ is a non-zero smooth mapping such that $k_1(t) \ell(t)+k_2(t) \beta(t)=0$ for all $t \in I$, $(\widetilde{k}_1, \widetilde{k}_2)$ is also a non-zero smooth mapping and 
\begin{align*}
\widetilde{k}_1(u)\widetilde{\ell}(u)+\widetilde{k}_2(u)\widetilde{\beta}(u)&=(k_1\circ t)(u)(\ell \circ t(u))t'(u)+(k_2 \circ t)(u)(\beta \circ t(u))t'(u)\\
&=t'(u)((k_1(t(u))\ell(t(u))+k_2(t(u))\beta(t(u)))=0.
\end{align*}
By a direct calculation, 
\begin{align*}
\dot{\widetilde{k}}_1(u) \widetilde{k}_2(u)-\widetilde{k}_1(u)\dot{\widetilde{k}}_2(u)&=(\dot k_1(t(u))t'(u))(k_2\circ t)(u)-(k_1\circ t)(u)(\dot k_2(t(u))t'(u))\\
&=t'(u)(\dot{k}_1(t(u)) k_2(t(u))-k_1(t(u)) \dot{k}_2(t(u))).
\end{align*}
Therefore, $u_0 \in \widetilde{I}$ is a vertex of $\widetilde{\gamma}$ if and only if $t(u_0) \in I$ is a vertex of $\gamma$.
\enD
Let $(\gamma,\nu):I \to \R^2 \times S^1$ be a Legendre curve with curvature $(\ell,\beta)$ and 
there exists a non-zero smooth mapping $(k_1,k_2):I \to \R^2 \setminus \{0\}$ such that $k_1(t) \ell(t)+k_2(t) \beta(t)=0$ for all $t \in I$.
We define $(\gamma^c,\nu^c):I \to \R^2 \times S^1$ by  
$$
\gamma^c(t)=\left(\int \ell(t) dt, \int \beta(t) dt\right), \ \nu^c(t)=(\cos \theta(t),\sin \theta(t)),
$$
where
$\cos \theta(t)={k_1(t)}/{\sqrt{k_1^2(t)+k_2^2(t)}}, \sin \theta(t)={k_2(t)}/{\sqrt{k_1^2(t)+k_2^2(t)}}.$
We call $\gamma^c$ a {\it curvature curve} of the Legendre curve $(\gamma,\nu)$.
\begin{proposition}\label{curvature-function}
Under the above notations, $(\gamma^c,\nu^c):I \to \R^2 \times S^1$ is a Legendre curve with the curvature $(\ell^c(t),\beta^c(t))=(\dot{\theta}(t), -\ell(t)\sin \theta(t)+\beta(t)\cos \theta(t))$. 
\end{proposition}
\demo
By a direct calculation, 
$\dot{\gamma}^c(t) \cdot \nu^c(t)=\ell(t) \cos \theta(t)+\beta(t) \sin \theta(t)=0$ for all $t \in I$. 
Hence, $(\gamma^c,\nu^c)$ is a Legendre curve. 
Since $\mu^c(t)=(-\sin \theta(t),\cos \theta(t))$, the curvature is given by 
$
\ell^c(t) =\dot{\theta}(t), \ \beta^c(t)=-\ell(t)\sin \theta(t)+\beta(t)\cos \theta(t).
$
\enD
By differentiating $\cos \theta(t)={k_1(t)}/{\sqrt{k_1^2(t)+k_2^2(t)}}, \sin \theta(t)={k_2(t)}/{\sqrt{k_1^2(t)+k_2^2(t)}}$, we have 
\begin{align*}
-\dot{\theta}(t)\sin \theta(t) &= \frac{k_2(t)(\dot{k}_1(t)k_2(t)-k_1(t)\dot{k}_2(t))}{(k_1^2(t)+k_2^2(t))^{\frac{3}{2}}}, \\
\dot{\theta}(t)\cos \theta(t) &=-\frac{k_1(t)(\dot{k}_1(t)k_2(t)-k_1(t)\dot{k}_2(t))}{(k_1^2(t)+k_2^2(t))^{\frac{3}{2}}}.
\end{align*}
It follows that $\dot{\theta}(t)=-(\dot{k}_1(t)k_2(t)-k_1(t)\dot{k}_2(t))/(k_1^2(t)+k_2^2(t))$.
By Proposition \ref{curvature-function}, we have the following corollary.
\begin{corollary}
Under the same notations as in Proposition \ref{curvature-function}, 
$t_0$ is a vertex of $\gamma$ if and only if $t_0\in I$ is an inflection point of the curvature curve $\gamma^c$.
\end{corollary}

As the corresponding result of Proposition \ref{condition2}, we have the following.
\begin{proposition}\label{condition3}
Let $(\gamma, \nu) : I \to \mathbb{R}^2 \times S^1$ be a Legendre curve with curvature $(\ell, \beta)$ and $\ell=\alpha \beta$, where $\alpha:I \to \R$ is a smooth function. 
\par
$(1)$ In the case of $n+1=k$. 
If $\gamma$ has four points of type $(n,m)=(n,2n+1)$, then $\gamma$ has at least four vertices. 
\par
$(2)$ In the case of $n+2 \le k$. 
If $\gamma$ has two points of type $(n,m)$, then $\gamma$ has at least four vertices. 
\end{proposition}
\demo
Suppose that a Legendre curve $(\gamma,\nu):I \to \R^2 \times S^1$ at $t_0 \in I$ is $\mathcal{R}$-equivalent to $(\widetilde{\gamma},\widetilde{\nu}):(\R,0) \to \R^2 \times S^1$,
\begin{align*}
\widetilde{\gamma}(t) =(\pm t^n,t^mf(t)),\
\widetilde{\nu}(t) =\frac{(-mt^kf(t)-t^{k+1}\dot{f}(t), \pm n)}{\sqrt{(mt^kf(t)+t^{k+1}\dot{f}(t))^2+n^2}},
\end{align*}
where $f:(\R,0) \to \R$ is a smooth function germ and the curvature $(\widetilde{\ell},\widetilde{\beta})$ (see Example \ref{nm-type}). 
If $n<k$, then there exists a unique function $\widetilde{\alpha}:(\R,0) \to \R$, 
$$
\widetilde{\alpha}(t)=\mp \frac{nt^{k-n}(mkf(t)+(m+k+1)t\dot{f}(t)+t^{2}\ddot{f}(t))}{((mt^kf(t)+t^{k+1}\dot{f}(t))^2+n^2)^{\frac{3}{2}}}
$$
such that $\widetilde{\ell}=\widetilde{\alpha} \widetilde{\beta}$. 
If $\gamma$ has a point of type $(n,m)$, then $\alpha$ at the point is $\mathcal{R}$-equivalent to $\widetilde{\alpha}$ at $0$. 
\par
$(1)$ If $n+1=k$, then $\widetilde{\alpha}$ is given by
$$
\widetilde{\alpha}(t)=\mp \frac{nt(m(n+1)f(t)+(m+n+2)t\dot{f}(t)+t^{2}\ddot{f}(t))}{((mt^{n+1}f(t)+t^{n+2}\dot{f}(t))^2+n^2)^{\frac{3}{2}}}.
$$
Since $\widetilde{\alpha}(0)=0$, we have $\alpha(t_0)=0$. 
If there exist two points such that $\alpha=0$, it is sufficient to show that there is at least a point such that $\dot\alpha=0$ between two adjacent $\alpha=0$. 
By the same argument of Proposition \ref{condition1}, $\gamma$ has at least four vertices. 
\par
$(2)$ If $n+2 \le k$, then $\widetilde{\alpha}$ and $\dot{\widetilde{\alpha}}$ are given by
\begin{align*}
\widetilde\alpha(t)=\mp &\frac{nt^{k-n}(mkf(t)+(m+k+1)t\dot{f}(t)+t^{2}\ddot{f}(t))}{((mt^kf(t)+t^{k+1}\dot{f}(t))^2+n^2)^{\frac{3}{2}}}, \\
\dot{\widetilde\alpha}(t)=\mp &\frac{n t^{k-n-1}}{((mt^kf(t)+t^{k+1}\dot{f}(t))^2+n^2)^{\frac{5}{2}}}\\
&\Bigl(\bigl((mt^kf(t)+t^{k+1}\dot{f}(t))^2+n^2\bigr)\bigl((k-n)(mkf(t)+(m+k+1)t\dot{f}(t)+t^2\ddot{f}(t))\\
&\hspace{5mm}+(mk+m+k+1)t\dot{f}(t)+(m+k+3)t^2\ddot{f}(t)+t^3\dddot{f}(t)\bigr)\\
&\hspace{10mm}-3t^{2k}(mf(t)+t\dot{f}(t))(mkf(t)+(m+k+1)t\dot{f}(t)+t^2\ddot{f}(t))^2\Bigr).
\end{align*}
Since $\widetilde{\alpha}(0)=\dot{\widetilde{\alpha}}(0)=0$, we have $\alpha(t_0)=\dot{\alpha}(t_0)=0$. 
By the same argument of Proposition \ref{condition1} (1), $\gamma$ has at least four vertices.
\enD

In general case, we have the four vertex theorem of frontals by the proofs of Propositions \ref{condition2} and \ref{condition3}.
\begin{theorem}
Let $(\gamma, \nu) : I \to \mathbb{R}^2 \times S^1$ be a Legendre curve with curvature $(\ell, \beta)$. 
Suppose that there exists a non-zero smooth mapping $(k_1,k_2):I \to \R^2 \setminus \{0\}$ such that $k_1(t) \ell(t)+k_2(t) \beta(t)=0$ for all $t \in I$. 
If there exist four points of the following types $(n,m)$:
$$
{\rm (i)}\ n=k \ {\rm and}\ \dot{f}(0)=0, \quad {\rm (ii)}\ n \ge k+2, \quad {\rm (iii)}\ n+2 \le k,
$$
then $\gamma$ has at least four vertices.
\end{theorem}

\begin{figure}[htbp]
  \centering
      \includegraphics[scale=0.175]{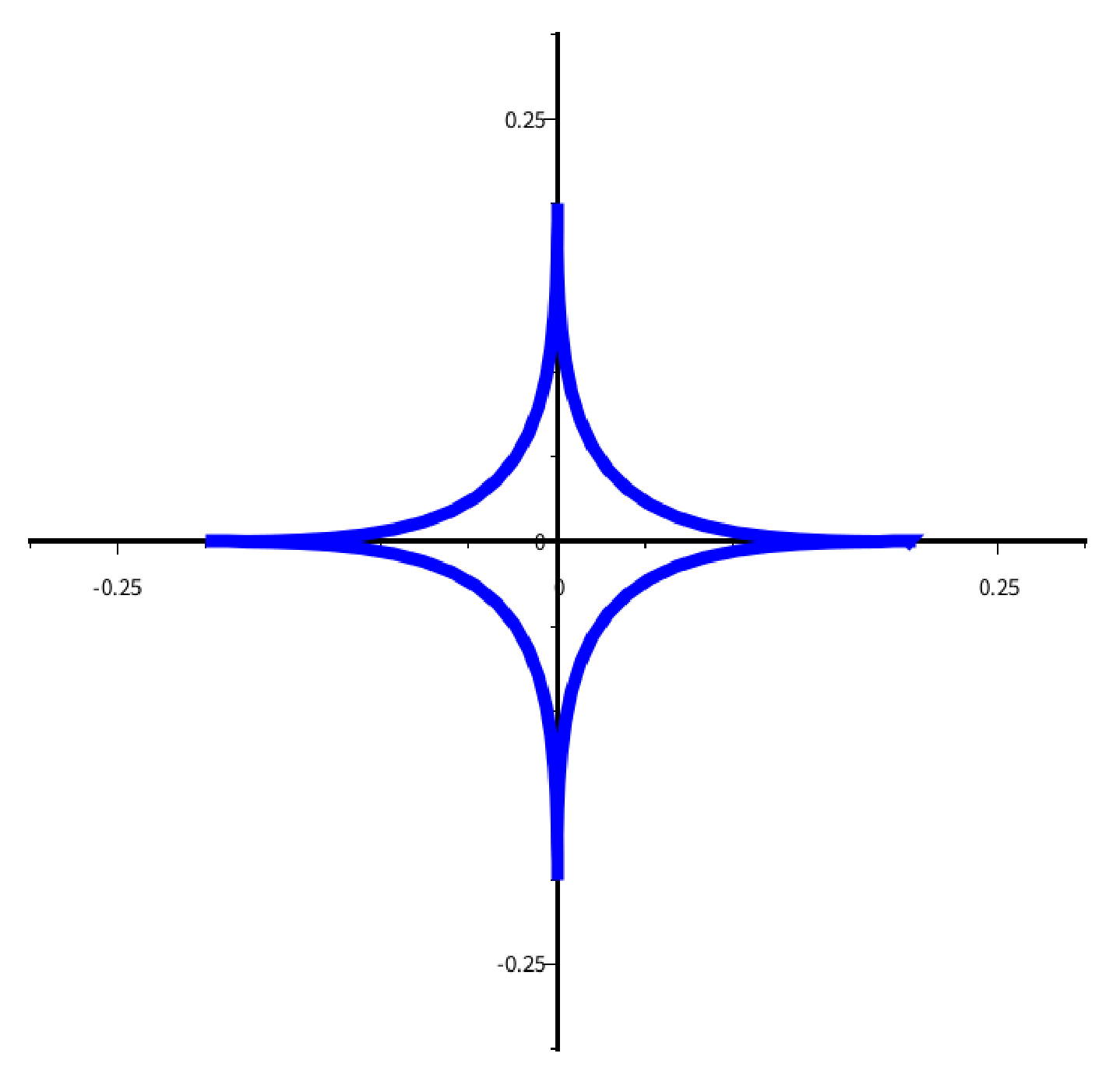}
      \includegraphics[scale=0.175]{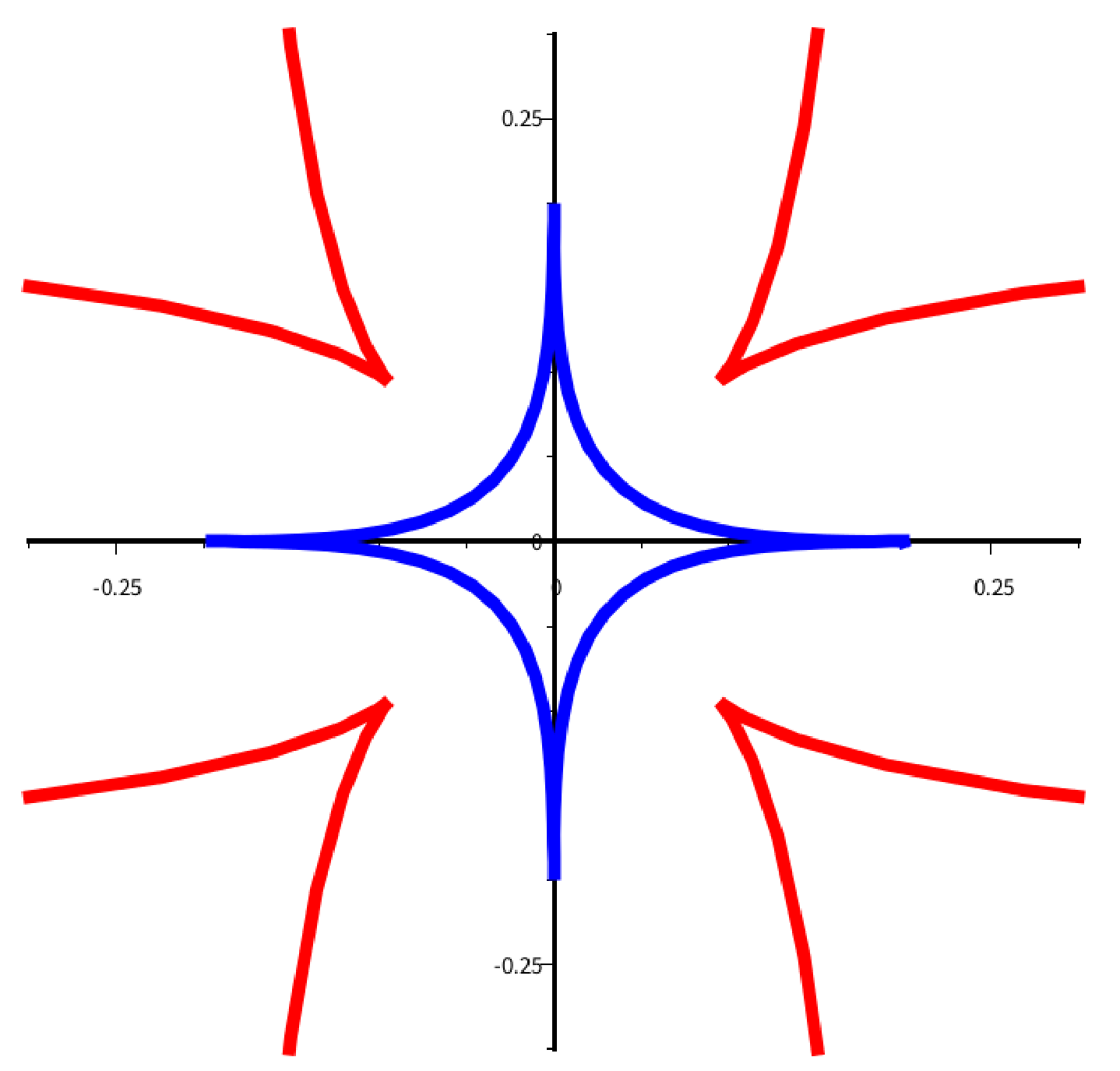}
  \caption{\centering{$\gamma(t)$ and $\mathcal{E}v(\gamma)(t)$}}
  \label{figure3}
\end{figure}

\begin{figure}[htbp]
  \centering
      \includegraphics[scale=0.185]{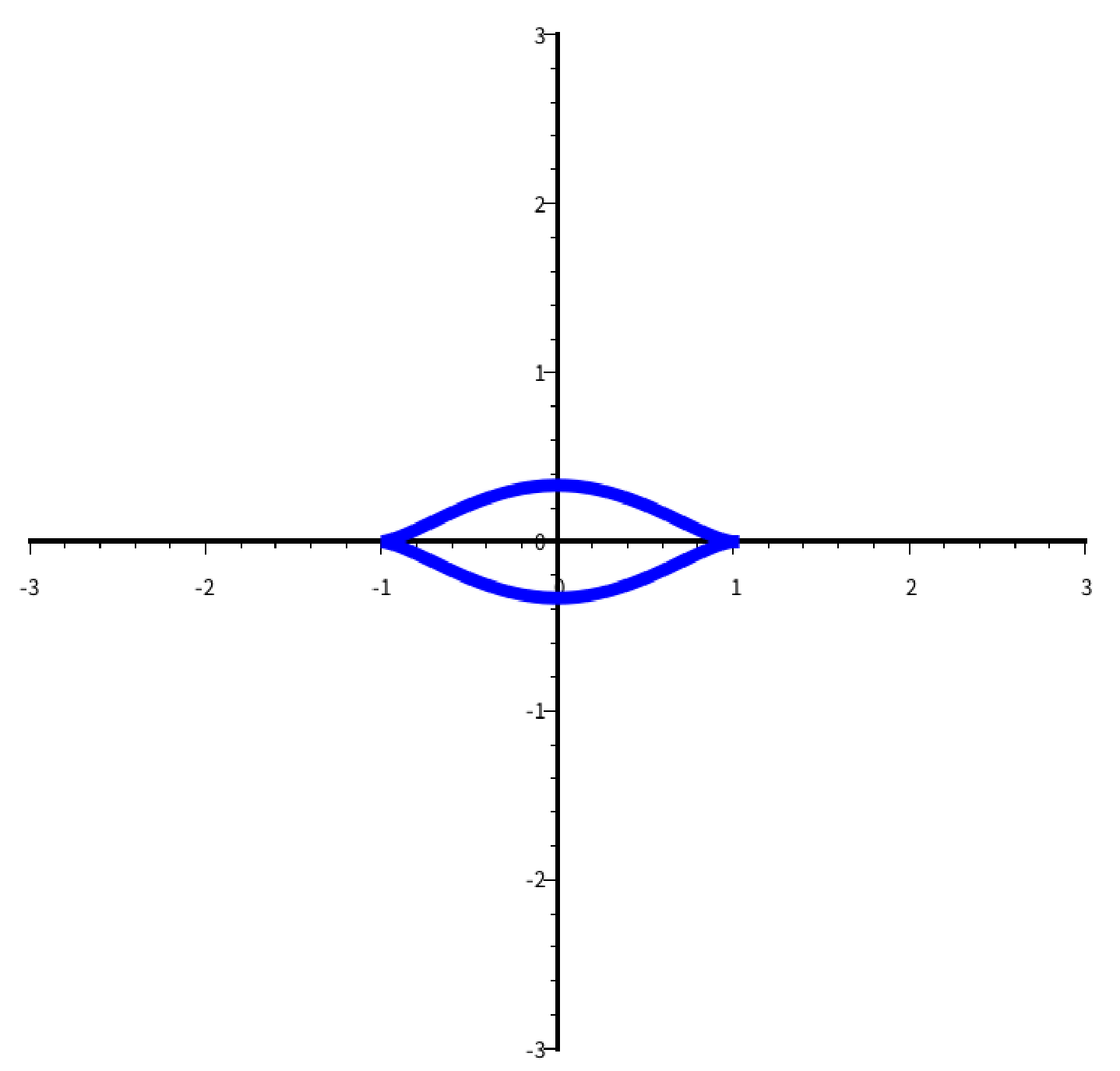}
      \includegraphics[scale=0.45]{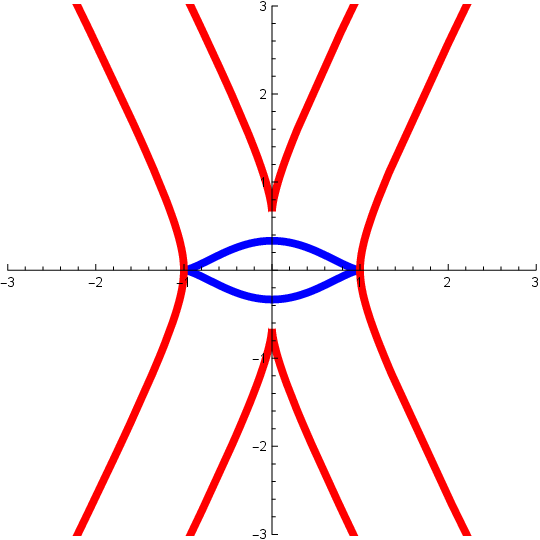}
  \caption{\centering{$\gamma(t)$ and $\mathcal{E}v(\gamma)(t)$}}
  \label{figure4}
\end{figure}

\begin{example}{\rm
Let $(\gamma,\nu):[0,2\pi] \to \R^2 \times S^1$ be 
\begin{align*}
 \gamma(t)=\left(\frac{1}{5}\cos^5 t, \frac{1}{5}\sin^5 t\right),\ \nu(t)=\frac{1}{\sqrt{\cos^6 t+\sin^6 t}}\left(-\sin^3 t, -\cos^3 t\right).
\end{align*}
By a direct calculation, $(\gamma,\nu)$ is a closed Legendre curve with the curvature
\begin{align*}
\ell(t)=-\frac{3\cos^2 t\sin^2 t}{\cos^6 t+\sin^6 t},  \
\beta(t)=-\cos t\sin t\sqrt{\cos^6 t+\sin^6 t}.
\end{align*}
If $k_1(t)=-(\cos^6 t+\sin^6 t)^{\frac{3}{2}}, \ k_2(t)=3\cos t\sin t$, then $(k_1,k_2)$ is a non-zero smooth function and $k_1(t) \ell(t)+k_2(t) \beta(t)=0$ for all $t \in I$. By a direct calculation, we have
\begin{align*}
\dot k_1(t)=9\cos 2t\cos t\sin t\sqrt{\cos^6 t+\sin^6 t}, \ \dot k_2(t)=3\cos 2t.
\end{align*}
Since
\begin{align*}
\dot k_1(t) k_2(t) -k_1(t)\dot k_2(t)=3\cos 2t\sqrt{\cos^6 t+\sin^6 t}(9\cos^2 t\sin^2 t+\cos^6 t+\sin^6 t),
\end{align*}
$\dot{k}_1(t) k_2(t) -k_1(t) \dot{k}_2(t)=0$ if and only if $\cos 2t=0$. 
It follows that $t=\pi/4, 3\pi/4, 5\pi/4, 7\pi/4$ are vertices of $\gamma$. 
Therefore, $\gamma$ has four vertices. 
Note that the singular points of $\gamma$ are $5/2$-cusps at $t=0, \pi/2, \pi, 3\pi/2$.
See Figure. \ref{figure3}}
\end{example}
\begin{example}{\rm
Let $(\gamma,\nu):[0,2\pi] \to \R^2 \times S^1$ be 
\begin{align*}
\gamma(t)=\left(\cos t, \frac{1}{3}\sin^3 t\right),\ \nu(t)=-\frac{1}{\sqrt{\cos^2 t\sin^2 t+1}}\left(\cos t\sin t,1\right).
\end{align*}
By a direct calculation, $(\gamma,\nu)$ is a closed Legendre curve with the curvature
\begin{align*}
\ell(t)=-\frac{\cos2t}{\cos^2 t\sin^2 t+1},\ 
\beta(t)=-\sin t\sqrt{\cos^2 t\sin^2 t+1}.
\end{align*} 
Note that $(\gamma, \nu)$ is a closed Legendre immersion, since $(\ell(t), \beta(t)) \neq (0, 0)$ for all $t \in [0,2\pi]$. 
By Theorem \ref{Convex-frontal}, $\gamma$ is not a convex front.
Moreover, $t=\pi/4, 3\pi/4, 5\pi/4, 7\pi/4$ are inflection points of $\gamma$.
If $k_1(t)=-\sin t(\cos^2 t\sin^2 t+1)^{\frac{3}{2}}, \ k_2(t)=\cos 2t$, then $(k_1,k_2)$ is a non-zero smooth function and $k_1(t) \ell(t)+k_2(t) \beta(t)=0$ for all $t\in I$. By a direct calculation, we have 
\begin{align*}
\dot k_1(t)=-\cos t(\cos^2 t\sin^2 t+1+3\sin^2 t\cos 2t)\sqrt{\cos^2 t\sin^2 t+1}, \dot k_2(t)=-2\sin 2t.
\end{align*}
Since
\begin{align*}
&\dot k_1(t) k_2(t) -k_1(t)\dot k_2(t)=-\cos t\sqrt{\cos^2 t\sin^2 t+1}\\
&\hspace{30mm}(\cos 2t(\cos^2 t\sin^2 t+1+3\sin^2 t\cos 2t)+4\sin^2 t(\cos^2 t\sin^2 t+1)),
\end{align*}
$\dot{k}_1(t) k_2(t) -k_1(t)\dot{k}_2(t)=0$ if and only if $\cos t=0$. 
It follows that $\gamma$ has only two vertices. 
Therefore, it does not satisfy on the four vertex theorem. 
Note that the singular points of $\gamma$ are $3/2$-cusps at $t=0, \pi$. See Figure. \ref{figure4}
}
\end{example}
\begin{theorem}\label{four-vertex-theorem-front-inflection}
Let $(\gamma,\nu):I \to \R^2 \times S^1$ be a closed Legendre immersion with curvature $(\ell,\beta)$. 
If $\gamma$ is a simple convex front with at least two singular and two inflection points, then $\gamma$ has at least four vertices.
\end{theorem}
\demo
Since $(\gamma,\nu)$ is a Legendre immersion, $(\ell(t),\beta(t)) \not= (0,0)$ for all $t \in I$.
If $k_1(t)=-\beta(t), k_2(t)=\ell(t)$, then $(k_1,k_2)$ is a non-zero smooth function and $k_1(t) \ell(t)+k_2(t)\beta(t)=0$ for all $t \in I$. 
Then we have $\dot{k}_1(t)k_2(t)-k_1(t)\dot{k}_2(t)=-\dot{\beta}(t)\ell(t)+\beta(t)\dot{\ell}(t)$. 
If $\ell(t_0)=0$ and $\dot{\ell}(t_0) \not=0$, then $\beta(t_0) \not=0$ and hence $t_0$ is an ordinary inflection point (cf. \cite{Bruce-Giblin,Gibson}). 
However, it contradicts the fact that $\gamma$ is convex.
Hence, if $\ell(t_0)=0$, then $\dot{\ell}(t_0)=0$. 
It follows that $t_0$ is a vertex of $\gamma$. 
Moreover, if $\beta(t_0)=0$ and $\dot{\beta}(t_0) \not=0$, then $\ell(t_0) \not=0$ and hence $\gamma$ has the point $t_0$ of type $(2,3)$. 
That is, an ordinary cusp ($3/2$-cusp) at $t_0$ (cf. \cite{Fukunaga-Takahashi-2014}).
However, it also contradicts the fact that $\gamma$ is convex.
Hence, if $\beta(t_0)=0$, then $\dot{\beta}(t_0)=0$. 
It follows that $t_0$ is a vertex of $\gamma$. 
Therefore, if $\gamma$ is a simple convex front with at least two singular and two inflection points, then $\gamma$ has at least four vertices.
\enD

We give a problem asking how about simple closed frontals with inflection points.
That is, does the four vertex theorem hold for a simple closed (convex) frontal with inflection points?

Nozomi Nakatsuyama, 
\\
Muroran Institute of Technology, Muroran 050-8585, Japan,
\\
E-mail address: 23043042@muroran-it.ac.jp
\\
\\
Masatomo Takahashi, 
\\
Muroran Institute of Technology, Muroran 050-8585, Japan,
\\
E-mail address: masatomo@muroran-it.ac.jp

\end{document}